\documentclass[journal,12pt,onecolumn,draftclsnofoot,]{IEEEtran}

\pdfoutput=1
\IEEEoverridecommandlockouts                              

\usepackage{graphicx}
\usepackage{amsmath,bm}
\usepackage{amssymb}
\usepackage{algorithm}
\usepackage{algpseudocode}
\usepackage{cite}
\usepackage{hyperref}

\usepackage[caption=false]{subfig}
\hypersetup{hidelinks}
\allowdisplaybreaks

\newtheorem{Remark}{Remark}

\newtheorem{Definition}{Definition}

\newenvironment{Proof}{\noindent{\em Proof:\/}}{\hfill \IEEEQED\par}
\newtheorem{Theorem}{Theorem}
\newtheorem{Lemma}{Lemma}

\newtheorem{Assumption}{Assumption}

\algdef{SE}[DOWHILE]{Do}{doWhile}{\algorithmicdo}[1]{\algorithmicwhile\ #1}%

\newcommand{\mathactivatecomma}{%
  \begingroup\lccode`~=`\,
  \lowercase{\endgroup\edef~}{\mathchar\the\mathcode`\,\penalty0 }}

\algnewcommand{\Initialize}[1]{%
  \State \textbf{Initialize: $j\in \mathcal{V}^i, i \in \mathcal{N}$}
  \Statex \hspace*{\algorithmicindent}\parbox[t]{.8\linewidth}{\raggedright #1}
}

\algnewcommand{\Iteration}[1]{%
  \State \textbf{Iteration $(t\geq 0)$: $j\in \mathcal{V}^i, i \in \mathcal{N}$}
  \Statex \hspace*{\algorithmicindent}\parbox[t]{.8\linewidth}{\raggedright #1}
}

\algnewcommand{\Output}[1]{%
  \State \textbf{Output: $j\in \mathcal{V}^i, i \in \mathcal{N}$}
  \Statex \hspace*{\algorithmicindent}\parbox[t]{.8\linewidth}{\raggedright #1}
}

\title{\LARGE \bf
Gradient-Free Nash Equilibrium Seeking in N-Cluster Games with Uncoordinated Constant Step-Sizes
}

\author{Yipeng Pang and Guoqiang Hu
\thanks{This research was supported by Singapore Ministry of Education Academic Research Fund Tier 1 RG180/17(2017-T1-002-158).}
\thanks{Y. Pang and G. Hu are with the School of Electrical and Electronic Engineering, Nanyang
Technological University, 639798, Singapore
        {\tt\small ypang005@e.ntu.edu.sg, gqhu@ntu.edu.sg}.}%
}

\begin{document}

\bstctlcite{IEEEexample:BSTcontrol}

\maketitle
\thispagestyle{empty}
\pagestyle{empty}

\begin{abstract}
This work investigates a problem of simultaneous global cost minimization and Nash equilibrium seeking, which commonly exists in $N$-cluster non-cooperative games. Specifically, the agents in the same cluster collaborate to minimize a global cost function, being a summation of their individual cost functions, and jointly play a non-cooperative game with other clusters as players.
For the problem settings, we suppose that the explicit analytical expressions of the agents' local cost functions are unknown, but the function values can be measured. We propose a gradient-free Nash equilibrium seeking algorithm by a synthesis of Gaussian smoothing techniques and gradient tracking. Furthermore, instead of using the uniform coordinated step-size, we allow the agents across different clusters to choose different constant step-sizes. When the largest step-size is sufficiently small, we prove a linear convergence of the agents' actions to a neighborhood of the unique Nash equilibrium under a strongly monotone game mapping condition, with the error gap being propotional to the largest step-size and the smoothing parameter. 
The performance of the proposed algorithm is validated by numerical simulations.
\end{abstract}


\section{Introduction}
The research on cooperation and competition across multiple interacting agents has been extensively studied in recent years, especially on distributed optimization and Nash equilibrium (NE) seeking in non-cooperative games. Specifically, distributed optimization deals with a cooperative minimization problem among a network of agents. 
On the other hand, NE seeking in non-cooperative games is concerned with a number of agents (also known as players), who are self-interested to minimize their individual cost given the other agents' decisions.

To simultaneously model the cooperative and competitive behaviors in networked systems, an $N$-cluster game is formulated. This game is essentially a non-cooperative game played among $N$ interacting clusters with each cluster being a virtual player. In each cluster, there are a group of agents who collaboratively minimize a cluster-level cost function given by a summation of their individual local cost functions. With these features, the $N$-cluster game naturally accommodates both collaboration and competition in a unified framework, which motivates us to study and propose solutions to find its NE. In this paper, we consider such an $N$-cluster non-cooperative game. Moreover, we further suppose that the explicit analytical expressions of the agents' local cost functions are unknown, but the function values can be measured.

A substantial works on NE seeking algorithms for non-cooperative games have been reported in the recent literature, including \cite{Ye2018b,Lu2019,Yi2019,DePersis2020,Zhang2019}, to list a few.
The focus of the aforementioned works is mainly on the competitive nature in the non-cooperative games. Different from that, the works in \cite{Gharesifard2013,Lou2016} considered two sub-networks zero-sum games, where each subnetwork owns an opposing cost function to be cooperatively minimized by the agents in the corresponding subnetwork. Then, an extension of such problems to $N$ subnetworks was firstly formulated in \cite{Ye2018}, which is known as an $N$-cluster (or coalition) game. Then, this problem has received a high research interest recently, which includes \cite{Ye2017b,Ye2017c,Ye2019,Nian2021,Zeng2019,Sun2021a,Ye2020,Pang2022}.
Most of the above works focus on continuous-time based methods, such as gradient play \cite{Ye2017b,Ye2017c,Ye2019,Nian2021}, subgradient dynamics \cite{Zeng2019}, projected primal-dual dynamics \cite{Sun2021a}, and extremum-seeking techniques \cite{Ye2020}.
Our previous work in \cite{Pang2022} introduced a discrete-time NE seeking strategy based on a synthesis of gradient-free and gradient-tracking techniques. This paper revisits the $N$-cluster game, and aims to extend the results to uncoordinated step-sizes across different clusters.

\textit{Contributions:} As compared to the aforementioned relevant works, the contributions of this work can be summarized as follows. 1) In contrast to the problem setups in \cite{Ye2017b,Ye2017c,Ye2019,Nian2021,Zeng2019,Sun2021a}, we limit the agents on the access to the cost functions: no explicit analytical expressions but only the values of the local cost functions can be utilized in the update laws. In this case, no gradient information can be directly utilized in the design of the algorithm. Hence, gradient-free techniques are adopted in this work. 2) As compared to our prior work \cite{Pang2022}, we extend the gradient-tracking method to allow uncoordinated constant step-sizes across different clusters, which further reduces the coordination among players from different clusters. 3) The technical challenges of the convergence analysis brought by gradient tracking methods in games, and uncoordinated step-sizes are addressed in this work. For the convergence results: we obtain a linear convergence to a neighborhood of the unique NE with the error being proportional to the largest step-size and a smoothing parameter under appropriate settings.

\textit{Notations:}
We use $\mathbf{1}_m$ ($\mathbf{0}_m$) for an $m$-dimensional vector with all elements being 1 (0), and $I_m$ for an $m\times m$ identity matrix. For a vector $\pi$, we use $\text{diag}(\pi)$ to denote a diagonal matrix formed by the elements of $\pi$. For any two vectors $u, v$, their inner product is denoted by $\langle u, v \rangle$.
The transpose of $u$ is denoted by $u^\top$. 
Moreover, we use $\|u\|$ for its standard Euclidean norm, \textit{i.e.}, $\|u\| = \sqrt{\langle u, u\rangle}$.
For vector $u$, we use $[u]_i$ to denote its $i$-th entry. 
The transpose and spectral norm of a matrix $A$ are denoted by $A^\top$ and $\|A\|$, respectively. 
We use $\rho(A)$ to represent the spectral radius of a square matrix $A$.
The expectation operator is denoted by $\mathbb{E}[\cdot]$.

\section{Problem Statement}\label{sec:problem_formulation}

\subsection{Problem Formulation}

An $N$-cluster game is defined by $\Gamma(\mathcal{N},\{f^i\},\{\mathbb{R}^{n_i}\})$, where each cluster, indexed by $i\in\mathcal{N} \triangleq \{1, 2, \ldots, N\}$, consists of a group of agents, denoted by $\mathcal{V}^i \triangleq \{1,2,\ldots,n_i\}$. Denote $n \triangleq \sum_{i=1}^Nn_i$. These agents aim to minimize a cluster-level cost function $f^i:\mathbb{R}^n\to\mathbb{R}$, defined as 
\begin{equation*}
f^i(\mathbf{x}^i,\mathbf{x}^{-i}) \triangleq \frac1{n_i}\sum_{j=1}^{n_i}f^i_j(\mathbf{x}^i,\mathbf{x}^{-i}), \quad \forall i\in\mathcal{N},
\end{equation*}
where $f^i_j(\mathbf{x}^i,\mathbf{x}^{-i})$ is a local cost function of agent $j$ in cluster $i$, $\mathbf{x}^i \triangleq [x^{i\top}_1,\ldots,x^{i\top}_{n_i}]^\top \in \mathbb{R}^{n_i}$ is a collection of all agents' actions in cluster $i$ with $x^i_j\in\mathbb{R}$ being the action of agent $j$ in cluster $i$, and $\mathbf{x}^{-i}\in\mathbb{R}^{n-n_i}$ denotes a collection of all agents' actions except cluster $i$. Denote $\mathbf{x}\triangleq (\mathbf{x}^{i},\mathbf{x}^{-i}) = [\mathbf{x}^{1\top},\ldots,\mathbf{x}^{N\top}]^\top$.

\begin{Definition}
(NE of $N$-Cluster Games). A vector $\mathbf{x}^* \triangleq (\mathbf{x}^{i*},\mathbf{x}^{-i*})\in\mathbb{R}^n$ is said to be an NE of the $N$-cluster non-cooperative game $\Gamma(\mathcal{N},\{f^i\},\{\mathbb{R}^{n_i}\})$, if and only if
\begin{align*}
f^i(\mathbf{x}^{i*},\mathbf{x}^{-i*})\leq f^i(\mathbf{x}^i,\mathbf{x}^{-i*}), \quad \forall \mathbf{x}^i\in\mathbb{R}^n,\quad \forall i\in \mathcal{N}.
\end{align*}
\end{Definition} 

Within each cluster $i\in\mathcal{N}$, there is an underlying directed communication network, denoted by $\mathcal{G}_i(\mathcal{V}^i,\mathcal{E}^i)$ with an adjacency matrix $\mathcal{A}^i\triangleq[a^i_{jk}]\in\mathbb{R}^{n_i\times n_i}$. In particular, $a^i_{jk}>0$ if agent $j$ can directly pass information to agent $k$, and $a^i_{jk}=0$ otherwise. We suppose $a^i_{jj}>0, \forall j\in\mathcal{V}^i$. Regarding the communication network, the following standard assumption is supposed.

\begin{Assumption}\label{assumption_graph}
For $i\in \mathcal{N}$, the digraph $\mathcal{G}_i(\mathcal{V}^i,\mathcal{E}^i)$ is strongly connected. The adjacency matrix $\mathcal{A}^i$ is doubly-stochastic.
\end{Assumption}

Noting that $\sigma_{\mathcal{A}^i} \triangleq \|\mathcal{A}^i-\frac1{n_i}\mathbf{1}_{n_i}\mathbf{1}_{n_i}^\top\|<1$ \cite[Lemma~1]{Pu2020}, we define $\bar{\sigma} \triangleq \max_{i\in\mathcal{N}}\sigma_{\mathcal{A}^i}$ and $\varsigma \triangleq \max_{i\in\mathcal{N}}(1+\sigma_{\mathcal{A}^i}^2)/(1-\sigma_{\mathcal{A}^i}^2)$.

Moreover, we consider the scenario where the explicit analytical expressions of the agents' local cost functions are unknown, but the function values can be measured, similar to the settings in \cite{Ye2020,Pang2022,Pang2020a,Pang2020c}. Regarding the cost function, the following standard assumption is supposed.

\begin{Assumption}\label{assumption_local_f_lipschitz}
For each $j\in\mathcal{V}^i,i\in\mathcal{N}$, the local cost function $f^i_j(\mathbf{x}^i,\mathbf{x}^{-i})$ is convex in $\mathbf{x}^i$, and continuously differentiable in $\mathbf{x}$. The total gradient $\nabla f^i_j(\mathbf{x})$ is $L$-Lipschitz continuous in $\mathbf{x}$, \textit{i.e.}, for any $\mathbf{x},\mathbf{x}'\in\mathbb{R}^n$, $\|\nabla f^i_j(\mathbf{x}) - \nabla f^i_j(\mathbf{x}')\|\leq L\|\mathbf{x} - \mathbf{x}'\|$.
\end{Assumption}

The game mapping of $\Gamma(\mathcal{N},\{f^i\},\{\mathbb{R}^{n_i}\})$ is defined as
\begin{align*}
\Phi(\mathbf{x})\triangleq [\nabla_{\mathbf{x}^1}f^1(\mathbf{x})^\top,\ldots,\nabla_{\mathbf{x}^N}f^N(\mathbf{x})^\top]^\top.
\end{align*}
The following standard assumption on the game mapping $\Phi(\mathbf{x})$ is supposed.

\begin{Assumption}\label{assumption_game_mapping}
The game mapping $\Phi$ of game $\Gamma$ is strongly monotone with a constant $\chi>0$, \textit{i.e.}, for any $\mathbf{x},\mathbf{x}'\in\mathbb{R}^n$, we have $\langle \Phi(\mathbf{x})-\Phi(\mathbf{x}'), \mathbf{x}-\mathbf{x}' \rangle\geq \chi\|\mathbf{x}-\mathbf{x}'\|^2$.
\end{Assumption}
\begin{Remark}
It is known that under Assumptions~\ref{assumption_local_f_lipschitz} and \ref{assumption_game_mapping}, game $\Gamma$ admits a unique NE.
\end{Remark}

\subsection{Preliminaries}
This part presents some preliminary results on gradient-free techniques based on Gaussian smoothing \cite{Nesterov2017}. 

For $j\in\mathcal{V}^i$, $i\in\mathcal{N}$, a Gaussian-smoothed function of the local cost function $f^i_j(\mathbf{x})$ can be defined as
\begin{align}
  f^i_{j,\mu}(\mathbf{x}) \triangleq \mathbb{E}_{\zeta\sim\mathcal{N}(\mathbf{0}_n, I_n)}[f^i_j(\mathbf{x}+\mu\zeta)], \label{eq:cost_function_smoothed}
\end{align}
where $\zeta$ is generated from a Gaussian distribution $\mathcal{N}(\mathbf{0}_n, I_n)$, and $\mu \geq 0$ is a smoothing parameter. 

For each cluster $i\in\mathcal{N}$, the randomized gradient-free oracle of $f^i_j(\mathbf{x})$ for player $j$ with respect to agent $k$, $j,k\in\mathcal{V}^i$, $i\in\mathcal{N}$ at time step $t\geq0$ is developed as
\begin{align}
{g}^i_{jk}(\mathbf{x}_t) \triangleq \frac{f^i_j(\mathbf{x}_t+\mu\zeta^i_{j,t})-f^i_j(\mathbf{x}_t)}{\mu}[\zeta^i_{j,t}]^i_k,\label{grad_oracle}
\end{align}
where $[\zeta^i_{j,t}]^i_k$ denotes the $(\sum_{l=0}^in_l+k)$-th element of $\zeta^i_{j,t}$ with $n_0 = 0$, and $\zeta^i_{j,t}$ being player $j$'s own version of $\zeta$ at time step $t$, and $\mu>0$. The oracle \eqref{grad_oracle} is useful as it can correctly estimate the partial gradient of the Gaussian-smoothed cost function $\nabla_{x^i_k} f^i_{j,\mu}(\mathbf{x}_t)$.
The following results for $f^i_{j,\mu}(\mathbf{x})$ and ${g}^i_{jk}(\mathbf{x})$ can be readily obtained according to \cite{Nesterov2017}.

\begin{Lemma}\label{lemma:property_f_mu}
(see \cite{Nesterov2017}) Under Assumption~\ref{assumption_local_f_lipschitz}, for $\forall j,k\in\mathcal{V}^i, i\in\mathcal{N}$, the following properties hold.
\begin{enumerate}
\item The function $f^i_{j,\mu}(\mathbf{x})$ is convex in $\mathbf{x}^i$ and totally differentiable in $\mathbf{x}$.
\item The total gradient $\nabla f^i_{j,\mu}(\mathbf{x})$ is $L$-Lipschitz continuous in $\mathbf{x}$, \textit{i.e.}, $\forall \mathbf{x}, \mathbf{y}\in\mathbb{R}^n$, $\|\nabla f^i_{j,\mu}(\mathbf{x}) - \nabla f^i_{j,\mu}(\mathbf{y})\| \leq L\|\mathbf{x}-\mathbf{y}\|$; and satisfies that $\|\nabla f^i_{j,\mu}(\mathbf{x}) - \nabla f^i_j(\mathbf{x})\| \leq \frac12(n+3)^{\frac32}L\mu$.
\item The randomized gradient-free oracle ${g}^i_{jk}(\mathbf{x})$ satisfies that $\mathbb{E}[{g}^i_{jk}(\mathbf{x})] = \nabla_{x^i_k} f^i_{j,\mu}(\mathbf{x})$, and $\mathbb{E}[\|{g}^i_{jk}(\mathbf{x})\|^2]\leq 4(n+4)\|\nabla f^i_{j,\mu}(\mathbf{x})\|^2+3(n+4)^3\mu^2L^2$.
\end{enumerate}
\end{Lemma}

We define a Gaussian-smoothed game associated with the $N$-cluster game $\Gamma$, denoted by $\Gamma_\mu(\mathcal{N},\{f^i_\mu\},\{\mathbb{R}^{n_i}\})$, having the same set of clusters and action sets as game $\Gamma$, but the cost function is given by
\begin{align*}
f^i_{\mu}(\mathbf{x}^i,\mathbf{x}^{-i}) \triangleq \frac1{n_i}\sum_{j=1}^{n_i}f^i_{j,\mu}(\mathbf{x}^i,\mathbf{x}^{-i}), \quad \forall i\in\mathcal{N},
\end{align*}
where $f^i_{j,\mu}$ is a Gaussian-smoothed function of $f^i_j$ defined in \eqref{eq:cost_function_smoothed}.
Similar to the game mapping of $\Gamma$, we define the game mapping of $\Gamma_\mu$ by $\Phi_\mu(\mathbf{x}) \triangleq [\nabla_{\mathbf{x}^1}f^1_\mu(\mathbf{x})^\top,\ldots,\nabla_{\mathbf{x}^N}f^N_\mu(\mathbf{x})^\top]^\top$.
The following lemma shows the strong monotonicity condition of $\Phi_\mu(\mathbf{x})$, and quantifies the distance between the NE of the smoothed game $\Gamma_\mu$ and the NE of the original game $\Gamma$ in terms of the smoothing parameter $\mu$.

\begin{Lemma}\label{lemma:NE_gap}
(see \cite[Lemma~1]{Pang2022}) Under Assumptions~\ref{assumption_local_f_lipschitz} and \ref{assumption_game_mapping}, for $\forall\mu\geq 0$, the smoothed game $\Gamma_\mu(\mathcal{N},\{f^i_\mu\},\{\mathbb{R}^{n_i}\})$ holds that
\begin{enumerate}
\item The game mapping $\Phi_\mu(\mathbf{x})$ is $\chi$-strongly monotone.
\item The smoothed game $\Gamma_\mu$ admits a unique NE (denoted by $\mathbf{x}_\mu^*$) satisfying that
\begin{align*}
\|\mathbf{x}_\mu^* - \mathbf{x}^*\| \leq \frac{n(n+3)^{\frac32}L\gamma}{2(1-\sqrt{1-\gamma\chi})}\mu,
\end{align*}
where $\mathbf{x}^*$ is the unique NE of the original game $\Gamma$, and $\gamma\in(0,\frac{\chi}{n^2L^2}]$ is a constant.
\end{enumerate}
\end{Lemma}

It follows from Lemma~\ref{lemma:NE_gap} that $\mathbf{x}^*_\mu$ is the unique NE of the smoothed game $\Gamma_\mu(\mathcal{N},\{f^i_\mu\},\{\mathbb{R}^{n_i}\})$, and hence holds that $\Phi_\mu(\mathbf{x}^*_\mu) = \mathbf{0}_n$. We define $G\triangleq \max_{j\in\mathcal{V}^i,i\in\mathcal{N}}\|\nabla f^i_{j,\mu}(\mathbf{x}^*_\mu)\|$.

\section{NE Seeking Algorithm for N-Cluster Games}\label{sec:distr_opt}

\subsection{Algorithm}

In this part, we present an NE seeking strategy for the $N$-Cluster Game. Specifically, each agent $j\in\mathcal{V}^i$, $i\in\mathcal{N}$ needs to maintain its own action variable $x^i_j$, and gradient tracker variables $\varphi^i_{jk}$ for $\forall k\in\mathcal{V}^i$. 
Let $x^i_{j,t}, \varphi^i_{jk,t}$ denote the values of these variables at time-step $t$. The update laws for each agent $j\in\mathcal{V}^i$, $i\in\mathcal{N}$ are designed as
\begin{subequations}\label{eq:algorithm}
\begin{align}
y^i_{jk,t+1} &= \sum_{l=1}^{n_i}a^i_{jl}y^i_{lk,t}-\alpha^i \varphi^i_{jk,t},\label{eq:update_y}\\
x^i_{j,t+1} &= y^i_{jj,t+1},\label{eq:update_x}\\
\varphi^i_{jk,t+1}& = \sum_{l=1}^{n_i}a^i_{jl} \varphi^i_{lk,t} + {g}^i_{jk}(\mathbf{x}_{t+1}) - {g}^i_{jk}(\mathbf{x}_t),\label{eq:update_grad_track}
\end{align}
\end{subequations}
with arbitrary $x^i_{j,0},y^i_{jk,0}\in\mathbb{R}$ and $\varphi^i_{jk,0} = {g}^i_{jk}(\mathbf{x}_0)$,
where ${g}^i_{jk}(\mathbf{x}_t)$ is the gradient estimator given by \eqref{grad_oracle}.
and $\alpha^i > 0$ is a constant step-size sequence adopted by all agents in cluster $i\in\mathcal{N}$. Denote the largest step-size by ${\alpha}_{\max}\triangleq\max_{i\in\mathcal{N}}\alpha^i$ and the average of all step-sizes by $\bar{\alpha} \triangleq \frac1n\sum_{i\in\mathcal{N}}n_i\alpha^i$.
Define the heterogeneity of the step-size as the following ratio, $\epsilon_\alpha\triangleq \|\bm{\alpha}-\bar{\bm{\alpha}}\|/\|\bar{\bm{\alpha}}\|$, where $\bm{\alpha} \triangleq [\alpha^1\mathbf{1}^\top_{n_1},\ldots,\alpha^N\mathbf{1}^\top_{n_N}]^\top$ and $\bar{\bm{\alpha}} \triangleq \bar{\alpha}\mathbf{1}_n$. 


\subsection{Main Results}

This part presents the main results of the proposed algorithm, as stated in the following theorem. Detailed proof is given in Sec.~\ref{subsec:proof_of_theorem}.
\begin{Theorem}\label{theorem:optimality}
Suppose Assumptions~\ref{assumption_graph}, \ref{assumption_local_f_lipschitz} and \ref{assumption_game_mapping} hold. Generate the auxiliary variables $\{y^i_{jk,t}\}_{t\geq0}$, the agent's action $\{x^i_{j,t}\}_{t\geq0}$ and gradient tracker $\{\varphi^i_{jk,t}\}_{t\geq0}$ by \eqref{eq:algorithm} with the uncoordinated constant step-size $\alpha^i$ satisfying $\epsilon_{\alpha} < \frac{\chi}{2\sqrt{n}L}$ and
\begin{align*}
0<{\alpha}_{\max}<\min\bigg\{\alpha_1,\alpha_2,\alpha_3,\frac1{\chi-2\sqrt{n}L\epsilon_{\alpha}},1\bigg\},
\end{align*}
where $\alpha_1$, $\alpha_2$ and $\alpha_3$ are defined in Sec.~\ref{subsec:proof_of_theorem}.
Then, all players' decisions $\mathbf{x}_t$ linearly converges to a neighborhood of the unique NE $\mathbf{x}^*$, and
\begin{align*}
\limsup_{t\to\infty}\mathbb{E}[\|\mathbf{x}_t - \mathbf{x}^*\|^2]\leq\mathcal{O}({\alpha}_{\max})+\mathcal{O}(\mu).
\end{align*}
\end{Theorem}

\begin{Remark}
Theorem~\ref{theorem:optimality} characterizes the convergence performance of the proposed algorithm. It shows that the agents' actions converge to a neighborhood of the NE linearly with the error bounded by two terms: one is proportional to the largest step-size, and the other is proportional to the smoothing parameter due to the gradient estimation.
\end{Remark}

\section{Convergence Analysis} \label{sec:conv_analysis}
Let $\mathcal{H}_t$ denote the $\sigma$-field generated by the entire history of the random variables from time-step 0 to $t-1$.
We introduce the following notations. Denote that $n_s\triangleq \sum_{i=1}^{N}n_i^2$ and $n_c\triangleq \sum_{i=1}^{N}n_i^3$. For $\forall k\in\mathcal{V}^i, i \in\mathcal{N}$,
$\mathbf{y}^i_{k,t} \triangleq [y^i_{1k,t}, \ldots, y^i_{n_ik,t}]^\top\in\mathbb{R}^{n_i},
\bar{y}^i_{k,t}\triangleq \frac1{n_i}\mathbf{1}_{n_i}^\top\mathbf{y}^i_{k,t}\in\mathbb{R},
\bar{\mathbf{y}}^i_t\triangleq [\bar{y}^i_{1,t},\ldots,\bar{y}^i_{n_i,t}]^\top\in\mathbb{R}^{n_i},
\bar{\mathbf{y}}_t\triangleq [\bar{\mathbf{y}}^{1\top}_t,\ldots,\bar{\mathbf{y}}^{N\top}_t]^\top\in\mathbb{R}^n,
\varphi^i_{k,t} \triangleq [\varphi^i_{1k,t},\ldots,\varphi^i_{n_ik,t}]^\top\in\mathbb{R}^{n_i}, 
\bar{\varphi}^i_{k,t} \triangleq \frac1{n_i}\mathbf{1}_{n_i}^\top\varphi^i_{k,t}\in\mathbb{R},
\bar{\bm{\varphi}}^i_t \triangleq [\bar{\varphi}^i_{1,t},\ldots,\bar{\varphi}^i_{n_i,t}]^\top\in\mathbb{R}^{n_i},
\mathbf{g}^i_{k} \triangleq [{g}^i_{1k},\ldots,{g}^i_{n_ik}]^\top\in\mathbb{R}^{n_i},
\bar{{g}}^i_k \triangleq \frac1{n_i}\mathbf{1}_{n_i}^\top\mathbf{g}^i_k\in\mathbb{R},
\nabla_{x^i_k} \mathbf{f}^i_\mu(\mathbf{x}) \triangleq [\nabla_{x^i_k} f^i_{1,\mu}(\mathbf{x}),\ldots,\nabla_{x^i_k} f^i_{n_i,\mu}(\mathbf{x})]^\top\in\mathbb{R}^{n_i}$.
Then, the update laws \eqref{eq:update_y} and \eqref{eq:update_grad_track} read:
\begin{subequations}
\begin{align}
\mathbf{y}^i_{k,t+1} &= \mathcal{A}^i\mathbf{y}^i_{k,t}-\alpha^i \varphi^i_{k,t},\label{eq:update_y_compact}\\
\varphi^i_{k,t+1} &= \mathcal{A}^i \varphi^i_{k,t} + \mathbf{g}^i_k(\mathbf{x}_{t+1}) - \mathbf{g}^i_k(\mathbf{x}_t).\label{eq:update_grad_track_compact}
\end{align}
\end{subequations}

Pre-multiplying both sides of \eqref{eq:update_y_compact} by $\frac1{n_i}\mathbf{1}_{n_i}^\top$ and augmenting the relation for $k\in\mathcal{V}^i$, we have
\begin{align}
\bar{\mathbf{y}}^i_{t+1} = \bar{\mathbf{y}}^i_t - \alpha^i\bar{\bm{\varphi}}^i_t, \label{eq:update_y_bar_compact}
\end{align}

The convergence analysis of the proposed algorithm will be conducted by: 1) constructing a linear system of three terms $\sum_{i=1}^N\sum_{k=1}^{n_i}\|\mathbf{y}^i_{k,t}-\mathbf{1}_{n_i}\bar{y}^i_{k,t}\|^2$, $\|\bar{\mathbf{y}}_\ell-\mathbf{x}^*_\mu\|^2$ and $\sum_{i=1}^N\sum_{k=1}^{n_i}\|\varphi^i_{k,t} - \mathbf{1}_{n_i}\bar{\varphi}^i_{k,t}\|^2$ in terms of their past iterations and some constants, 2) analyzing the convergence of the established linear system. 

\subsection{Auxiliary Analysis}

We first derive some results for the averaged gradient tracker $\bar{\varphi}^i_{k,t}$.
\begin{Lemma}\label{lemma:averaged_grad_tracker}
Under Assumptions~\ref{assumption_graph} and \ref{assumption_local_f_lipschitz}, the averaged gradient tracker $\bar{\varphi}^i_{k,t}, \forall k\in\mathcal{V}^i, i \in\mathcal{N}$ holds that
\begin{enumerate}
\item $\bar{\varphi}^i_{k,t} = \bar{{g}}^i_{k}(\mathbf{x}_t)$,
\item $\mathbb{E}[\bar{\varphi}^i_{k,t}|\mathcal{H}_t] = \nabla_{x^i_k} f^i_\mu(\mathbf{x}_t)$,
\item $\mathbb{E}[\|\bar{\varphi}^i_{k,t}\|^2|\mathcal{H}_t] \leq 12(n+4)L^2\sum_{i=1}^N\sum_{k=1}^{n_i}\|\mathbf{y}^i_{k,t}-\mathbf{1}_{n_i}\bar{y}^i_{k,t}\|^2+12(n+4)L^2\|\bar{\mathbf{y}}_t-\mathbf{x}^*_\mu\|^2+12(n+4)G^2+3(n+4)^3\mu^2L^2$.
\end{enumerate}
\end{Lemma}
\begin{Proof}
For 1), multiplying $\frac1{n_i}\mathbf{1}_{n_i}^\top$ from the left on both sides of \eqref{eq:update_grad_track_compact}, and noting that $\mathcal{A}^i$ is doubly stochastic, we have
\begin{align*}
\bar{\varphi}^i_{k,t+1} = \bar{\varphi}^i_{k,t} + \bar{{g}}^i_k(\mathbf{x}_{t+1}) - \bar{{g}}^i_k(\mathbf{x}_t).
\end{align*}
Recursively expanding the above relation and noting that $\varphi^i_{k,0} = \mathbf{g}^i_k(\mathbf{x}_0)$ completes the proof.

For 2), following the result of part 1) and Lemma~\ref{lemma:property_f_mu}-3), we obtain
\begin{align*}
&\mathbb{E}[\bar{\varphi}^i_{k,t}|\mathcal{H}_t] = \mathbb{E}[\bar{{g}}^i_k(\mathbf{x}_t)|\mathcal{H}_t]=\frac1{n_i}\mathbf{1}_{n_i}^\top\mathbb{E}[\mathbf{g}^i_k|\mathcal{H}_t] =\frac1{n_i}\mathbf{1}_{n_i}^\top\nabla_{x^i_k} \mathbf{f}^i_\mu=\nabla_{x^i_k} f^i_\mu(\mathbf{x}_t).
\end{align*}

For 3), it follows that
\begin{align*}
&\mathbb{E}[\|\bar{\varphi}^i_{k,t}\|^2|\mathcal{H}_t]=\mathbb{E}[\|\bar{{g}}^i_{k}(\mathbf{x}_t)\|^2|\mathcal{H}_t] = \frac1{n_i^2}\mathbb{E}[\|\mathbf{1}_{n_i}^\top\mathbf{g}^i_k(\mathbf{x}_t)\|^2|\mathcal{H}_t]\leq\frac1{n_i}\sum_{j=1}^{n_i}\mathbb{E}[\|{g}^i_{jk}(\mathbf{x}_t)\|^2|\mathcal{H}_t].
\end{align*}
On the other hand, it follows from Lemma~\ref{lemma:property_f_mu}-3) that
\begin{align}
\mathbb{E}[\|{g}^i_{jk}(\mathbf{x}_t)\|^2|\mathcal{H}_t]&\leq 4(n+4)\|\nabla f^i_{j,\mu}(\mathbf{x}_t)\|^2+3(n+4)^3\mu^2L^2\nonumber\\
&\leq 12(n+4)\|\nabla f^i_{j,\mu}(\mathbf{x}_t)-\nabla f^i_{j,\mu}(\bar{\mathbf{y}}_t)\|^2+12(n+4)\|\nabla f^i_{j,\mu}(\bar{\mathbf{y}}_t)-\nabla f^i_{j,\mu}(\mathbf{x}^*_\mu)\|^2\nonumber\\
&\quad+12(n+4)\|\nabla f^i_{j,\mu}(\mathbf{x}^*_\mu)\|^2+3(n+4)^3\mu^2L^2\nonumber\\
&\leq12(n+4)L^2\|\mathbf{x}_t-\bar{\mathbf{y}}_t\|^2+12(n+4)L^2\|\bar{\mathbf{y}}_t-\mathbf{x}^*_\mu\|^2\nonumber\\
&\quad+12(n+4)G^2+3(n+4)^3\mu^2L^2\nonumber\\
&\leq12(n+4)L^2\sum_{i=1}^N\sum_{k=1}^{n_i}\|\mathbf{y}^i_{k,t}-\mathbf{1}_{n_i}\bar{y}^i_{k,t}\|^2\nonumber\\
&\quad+12(n+4)(L^2\|\bar{\mathbf{y}}_t-\mathbf{x}^*_\mu\|^2+G^2)+3(n+4)^3\mu^2L^2,\label{eq:pi_i_jk_bound}
\end{align}
where $G\triangleq \max_{j\in\mathcal{V}^i,i\in\mathcal{N}}\|\nabla f^i_{j,\mu}(\mathbf{x}^*_\mu)\|$ and the last inequality follows from \eqref{eq:update_x} that
\begin{align}
\|\mathbf{x}_t-\bar{\mathbf{y}}_t\|^2&=\sum_{i=1}^N\sum_{k=1}^{n_i}\|y^i_{kk,t}-\bar{y}^i_{k,t}\|^2\leq\sum_{i=1}^N\sum_{k=1}^{n_i}\|\mathbf{y}^i_{k,t} -\mathbf{1}_{n_i}\bar{y}^i_{k,t}\|^2. \label{eq:x_minus_ybar}
\end{align}
The proof is completed by combining the above relations.
\end{Proof}

Then, we provide a bound on the stacked gradient tracker $\varphi^i_{k,t}$.
\begin{Lemma}\label{lemma:stacked_grad_tracker}
Under Assumptions~\ref{assumption_graph} and \ref{assumption_local_f_lipschitz}, the stacked gradient tracker $\{\varphi^i_{k,t}\}_{t\geq0}, \forall k\in\mathcal{V}^i, i \in\mathcal{N}$ holds that
\begin{align*}
\mathbb{E}[\|\varphi^i_{k,t}\|^2|\mathcal{H}_t]&\leq 2\mathbb{E}[\|\varphi^i_{k,t} - \mathbf{1}_{n_i}\bar{\varphi}^i_{k,t}\|^2|\mathcal{H}_t] + 24n_i^2(n+4)L^2\sum_{i=1}^N\sum_{k=1}^{n_i}\|\mathbf{y}^i_{k,t}-\mathbf{1}_{n_i}\bar{y}^i_{k,t}\|^2\\
&\quad+24n_i^2(n+4)(L^2\|\bar{\mathbf{y}}_t-\mathbf{x}^*_\mu\|^2+G^2)+6n_i^2(n+4)^3\mu^2L^2.
\end{align*}
\end{Lemma}
\begin{Proof}
It is noted that
\begin{align*}
\|\varphi^i_{k,t}\|^2 &\leq 2\|\varphi^i_{k,t} - \mathbf{1}_{n_i}\bar{\varphi}^i_{k,t}\|^2 + 2\|\mathbf{1}_{n_i}\bar{\varphi}^i_{k,t}\|^2=2\|\varphi^i_{k,t} - \mathbf{1}_{n_i}\bar{\varphi}^i_{k,t}\|^2 + 2n_i^2\|\bar{\varphi}^i_{k,t}\|^2.
\end{align*}
The proof is completed by taking the conditional expectation on $\mathcal{H}_t$ on both sides and substituting Lemma~\ref{lemma:averaged_grad_tracker}-3).
\end{Proof}

Now, we are ready to establish an inequality relation for the first term, $\sum_{i=1}^N\sum_{k=1}^{n_i}\|\mathbf{y}^i_{k,t}-\mathbf{1}_{n_i}\bar{y}^i_{k,t}\|^2$ in Lemma~\ref{lemma:conensus}.

\begin{Lemma}\label{lemma:conensus}
Under Assumptions~\ref{assumption_graph}, \ref{assumption_local_f_lipschitz} and \ref{assumption_game_mapping}, the total consensus error of the auxiliary variables $\sum_{i=1}^N\sum_{k=1}^{n_i}\|\mathbf{y}^i_{k,t}-\mathbf{1}_{n_i}\bar{y}^i_{k,t}\|^2$ satisfies that
\begin{align*}
&\sum_{i=1}^N\sum_{k=1}^{n_i}\mathbb{E}[\|\mathbf{y}^i_{k,t+1}-\mathbf{1}_{n_i}\bar{y}^i_{k,t+1}\|^2|\mathcal{H}_t] \leq\bigg(\frac{1+\bar{\sigma}^2}{2}+24(n+4)n_c\varsigma L^2\alpha_{\max}^2\bigg)\sum_{i=1}^N\sum_{k=1}^{n_i}\|\mathbf{y}^i_{k,t}-\mathbf{1}_{n_i}\bar{y}^i_{k,t}\|^2\\
&\quad+24(n+4)n_c\varsigma L^2\alpha_{\max}^2\|\bar{\mathbf{y}}_t-\mathbf{x}^*_\mu\|^2+2\varsigma\alpha_{\max}^2\sum_{i=1}^N\sum_{k=1}^{n_i}\mathbb{E}[\|\varphi^i_{k,t} - \mathbf{1}_{n_i}\bar{\varphi}^i_{k,t}\|^2|\mathcal{H}_t]\\
&\quad +24(n+4)n_c\varsigma G^2\alpha_{\max}^2+6(n+4)^3n_c\varsigma \mu^2L^2\alpha_{\max}^2.
\end{align*}
\end{Lemma}
\begin{Proof}
It follows from \eqref{eq:update_y_compact} that for $i\in\mathcal{N}$
\begin{align*}
\|\mathbf{y}^i_{k,t+1}-\mathbf{1}_{n_i}\bar{y}^i_{k,t+1}\|^2 &= \|\mathcal{A}^i\mathbf{y}^i_{k,t} -\alpha^i \varphi^i_{k,t}-\frac1{n_i}\mathbf{1}_{n_i}\mathbf{1}_{n_i}^\top (\mathcal{A}^i\mathbf{y}^i_{k,t} -\alpha^i \varphi^i_{k,t})\|^2\\
&\leq\|\mathcal{A}^i\mathbf{y}^i_{k,t}-\mathbf{1}_{n_i}\bar{y}^i_{k,t}\|^2 + \|\alpha^i(I_{n_i}-\frac1{n_i}\mathbf{1}_{n_i}\mathbf{1}_{n_i}^\top) \varphi^i_{k,t}\|^2 \\
&\quad-2\alpha^i\langle \mathcal{A}^i\mathbf{y}^i_{k,t}-\mathbf{1}_{n_i}\bar{y}^i_{k,t}, (I_{n_i}-\frac1{n_i}\mathbf{1}_{n_i}\mathbf{1}_{n_i}^\top)\varphi^i_{k,t}\rangle,
\end{align*}
Taking the conditional expectation on $\mathcal{H}_t$ and noting that $\|I_{n_i}-\frac1{n_i}\mathbf{1}_{n_i}\mathbf{1}_{n_i}^\top\| = 1$, we obtain
\begin{align*}
&\mathbb{E}[\|\mathbf{y}^i_{k,t+1}-\mathbf{1}_{n_i}\bar{y}^i_{k,t+1}\|^2|\mathcal{H}_t] \leq\sigma_{\mathcal{A}^i}^2\|\mathbf{y}^i_{k,t}-\mathbf{1}_{n_i}\bar{y}^i_{k,t}\|^2+ \alpha_{\max}^2\mathbb{E}[\|\varphi^i_{k,t}\|^2|\mathcal{H}_t]\\
&\quad\quad +\frac{1-\sigma_{\mathcal{A}^i}^2}{2\sigma_{\mathcal{A}^i}^2}\mathbb{E}[\|\mathcal{A}^i\mathbf{y}^i_{k,t}-\mathbf{1}_{n_i}\bar{y}^i_{k,t}\|^2|\mathcal{H}_t]+\frac{2\sigma_{\mathcal{A}^i}^2}{1-\sigma_{\mathcal{A}^i}^2}\alpha_{\max}^2\mathbb{E}[\|\varphi^i_{k,t}\|^2|\mathcal{H}_t]\\
&\quad\leq\frac{1+\sigma_{\mathcal{A}^i}^2}2\mathbb{E}[\|\mathbf{y}^i_{k,t}-\mathbf{1}_{n_i}\bar{y}^i_{k,t}\|^2|\mathcal{H}_t]+\frac{1+\sigma_{\mathcal{A}^i}^2}{1-\sigma_{\mathcal{A}^i}^2}\alpha_{\max}^2\mathbb{E}[\|\varphi^i_{k,t}\|^2|\mathcal{H}_t]\\
&\quad\leq\frac{1+\bar{\sigma}^2}{2}\|\mathbf{y}^i_{k,t}-\mathbf{1}_{n_i}\bar{y}^i_{k,t}\|^2+\varsigma\alpha_{\max}^2\mathbb{E}[\|\varphi^i_{k,t}\|^2|\mathcal{H}_t].
\end{align*}
Applying Lemma~\ref{lemma:stacked_grad_tracker} and summing over $k=1$ to $n_i$, $i=1$ to $N$ complete the proof.
\end{Proof}

Then, we proceed to build the inequality relation for the second term, $\|\bar{\mathbf{y}}_t-\mathbf{x}^*_\mu\|^2$ in Lemma~\ref{lemma:optimality_gap}.

\begin{Lemma}\label{lemma:optimality_gap}
Under Assumptions~\ref{assumption_graph}, \ref{assumption_local_f_lipschitz} and \ref{assumption_game_mapping}, the gap between the stacked averaged auxiliary variable and the NE of game $\Gamma_\mu$, $\|\bar{\mathbf{y}}_t-\mathbf{x}^*_\mu\|^2$ holds that
\begin{align*}
\mathbb{E}[\|\bar{\mathbf{y}}_{t+1}-\mathbf{x}^*\|^2|\mathcal{H}_t] &\leq(1-(\chi-2\sqrt{n}L\epsilon_{\alpha})\bar{\alpha}+12n(n+4)L^2\alpha_{\max}^2)\|\bar{\mathbf{y}}_t-\mathbf{x}^*_\mu\|^2\\
&\quad+\bigg(12n(n+4)L^2\alpha_{\max}^2+\frac{n^2L^2\alpha_{\max}}{\chi}\bigg)\sum_{i=1}^N\sum_{k=1}^{n_i}\|\mathbf{y}^i_{k,t}-\mathbf{1}_{n_i}\bar{y}^i_{k,t}\|^2\\
&\quad+12n(n+4)G^2\alpha_{\max}^2+3n(n+4)^3\mu^2L^2\alpha_{\max}^2.
\end{align*}
\end{Lemma}
\begin{Proof}
From \eqref{eq:update_y_bar_compact}, we know
\begin{subequations}\label{eq:optimality_gap_expand}
\begin{align}
\|\bar{\mathbf{y}}_{t+1}-\mathbf{x}^*_\mu\|^2 &= \sum_{i=1}^N\sum_{k=1}^{n_i}\|\bar{y}^i_{k,t} - x^{i*}_{k,\mu} - \alpha^i \bar{\varphi}^i_{k,t}\|^2\nonumber\\
&\leq\|\bar{\mathbf{y}}_t-\mathbf{x}^*_\mu\|^2 + \alpha_{\max}^2\sum_{i=1}^N\sum_{k=1}^{n_i}\|\bar{\varphi}^i_{k,t}\|^2\label{eq:optimality_gap_expand_1}\\
&\quad-2\sum_{i=1}^N\sum_{k=1}^{n_i}\alpha^i\langle\bar{y}^i_{k,t} - x^{i*}_{k,\mu}, \bar{\varphi}^i_{k,t}-\nabla_{x^i_k} f^i_\mu(\mathbf{x}_t)\rangle\label{eq:optimality_gap_expand_2}\\
&\quad-2\sum_{i=1}^N\sum_{k=1}^{n_i}\alpha^i\langle\bar{y}^i_{k,t} - x^{i*}_{k,\mu}, \nabla_{x^i_k} f^i_\mu(\mathbf{x}_t)-\nabla_{x^i_k} f^i_\mu(\bar{\mathbf{y}}_t)\rangle\label{eq:optimality_gap_expand_3}\\
&\quad-2\sum_{i=1}^N\sum_{k=1}^{n_i}\alpha^i\langle\bar{y}^i_{k,t} - x^{i*}_{k,\mu}, \nabla_{x^i_k} f^i_\mu(\bar{\mathbf{y}}_t) - \nabla_{x^i_k} f^i_\mu(\mathbf{x}^*_\mu)\rangle.\label{eq:optimality_gap_expand_4}
\end{align}
\end{subequations}
For the second term in \eqref{eq:optimality_gap_expand_1}, applying Lemma~\ref{lemma:averaged_grad_tracker}-(3)
\begin{align}
&\alpha_{\max}^2\sum_{i=1}^N\sum_{k=1}^{n_i}\mathbb{E}[\|\bar{\varphi}^i_{k,t}\|^2|\mathcal{H}_t]\leq\alpha_{\max}^2\sum_{i=1}^N\sum_{k=1}^{n_i}\bigg(12(n+4)L^2\sum_{i=1}^N\sum_{k=1}^{n_i}\|\mathbf{y}^i_{k,t}-\mathbf{1}_{n_i}\bar{y}^i_{k,t}\|^2\nonumber\\
&\quad\quad+12(n+4)(L^2\|\bar{\mathbf{y}}_t-\mathbf{x}^*_\mu\|^2+G^2)+3(n+4)^3\mu^2L^2\bigg)\nonumber\\
&\quad\leq12n(n+4)L^2\alpha_{\max}^2\sum_{i=1}^N\sum_{k=1}^{n_i}\|\mathbf{y}^i_{k,t}-\mathbf{1}_{n_i}\bar{y}^i_{k,t}\|^2+12n(n+4)L^2\alpha_{\max}^2\|\bar{\mathbf{y}}_t-\mathbf{x}^*_\mu\|^2\nonumber\\
&\quad\quad+12n(n+4)G^2\alpha_{\max}^2+3n(n+4)^3\mu^2L^2\alpha_{\max}^2. \label{eq:optimality_gap_expand_1_result}
\end{align}
For \eqref{eq:optimality_gap_expand_2},
\begin{align}
\mathbb{E}[-2\alpha^i\langle\bar{y}^i_{k,t} - x^{i*}_{k,\mu}, \bar{\varphi}^i_{k,t}-\nabla_{x^i_k} f^i_\mu(\mathbf{x}_t)\rangle|\mathcal{H}_t]=0. \label{eq:optimality_gap_expand_2_result}
\end{align}
For \eqref{eq:optimality_gap_expand_3}, 
\begin{align}
&-2\sum_{i=1}^N\sum_{k=1}^{n_i}\alpha^i\langle\bar{y}^i_{k,t} - x^{i*}_{k,\mu}, \nabla_{x^i_k} f^i_\mu(\mathbf{x}_t)-\nabla_{x^i_k} f^i_\mu(\bar{\mathbf{y}}_t)\rangle\nonumber\\
&\quad\leq2\alpha_{\max}\sum_{i=1}^N\sum_{k=1}^{n_i}\|\bar{y}^i_{k,t} - x^{i*}_{k,\mu}\|\|\nabla_{x^i_k} f^i_\mu(\mathbf{x}_t)-\nabla_{x^i_k} f^i_\mu(\bar{\mathbf{y}}_t)\|\nonumber\\
&\quad\leq2L\alpha_{\max}\sum_{i=1}^N\sum_{k=1}^{n_i}\|\bar{y}^i_{k,t} - x^{i*}_{k,\mu}\|\|\mathbf{x}_t-\bar{\mathbf{y}}_t\|\leq2\sqrt{n}L\alpha_{\max}\|\bar{\mathbf{y}}_t - \mathbf{x}^*_\mu\|\|\mathbf{x}_t-\bar{\mathbf{y}}_t\|\nonumber\\
&\quad\leq\chi\bar{\alpha}\|\bar{\mathbf{y}}_t - \mathbf{x}^*_\mu\|^2+\frac{nL^2\alpha_{\max}^2}{\chi\bar{\alpha}}\|\mathbf{x}_t-\bar{\mathbf{y}}_t\|^2\nonumber\\
&\quad\leq\chi\bar{\alpha}\|\bar{\mathbf{y}}_t - \mathbf{x}^*_\mu\|^2+\frac{nL^2\alpha_{\max}^2}{\chi\bar{\alpha}}\sum_{i=1}^N\sum_{k=1}^{n_i}\|\mathbf{y}^i_{k,t}-\mathbf{1}_{n_i}\bar{y}^i_{k,t}\|^2\nonumber\\
&\quad\leq\chi\bar{\alpha}\|\bar{\mathbf{y}}_t - \mathbf{x}^*_\mu\|^2+\frac{n^2L^2\alpha_{\max}}{\chi}\sum_{i=1}^N\sum_{k=1}^{n_i}\|\mathbf{y}^i_{k,t}-\mathbf{1}_{n_i}\bar{y}^i_{k,t}\|^2, \label{eq:optimality_gap_expand_3_result}
\end{align}
where the last inequality is due to $\frac{{\alpha}_{\max}}{\bar{\alpha}}<n$.
For \eqref{eq:optimality_gap_expand_4}, 
\begin{subequations}\label{eq:optimality_gap_expand_4_break}
\begin{align}
&-2\sum_{i=1}^N\sum_{k=1}^{n_i}\alpha^i\langle\bar{y}^i_{k,t} - x^{i*}_{k,\mu}, \nabla_{x^i_k} f^i_\mu(\bar{\mathbf{y}}_t) - \nabla_{x^i_k} f^i_\mu(\mathbf{x}^*_\mu)\rangle\nonumber\\
&\quad=-2\sum_{i=1}^N\sum_{k=1}^{n_i}(\alpha^i-\bar{\alpha})\langle\bar{y}^i_{k,t} - x^{i*}_{k,\mu}, \nabla_{x^i_k} f^i_\mu(\bar{\mathbf{y}}_t) - \nabla_{x^i_k} f^i_\mu(\mathbf{x}^*_\mu)\rangle \label{eq:optimality_gap_expand_4_break_1}\\
&\quad\quad-2\bar{\alpha}\sum_{i=1}^N\sum_{k=1}^{n_i}\langle\bar{y}^i_{k,t} - x^{i*}_{k,\mu}, \nabla_{x^i_k} f^i_\mu(\bar{\mathbf{y}}_t) - \nabla_{x^i_k} f^i_\mu(\mathbf{x}^*_\mu)\rangle. \label{eq:optimality_gap_expand_4_break_2}
\end{align}
\end{subequations}
For \eqref{eq:optimality_gap_expand_4_break_1},
\begin{align}
&-2\sum_{i=1}^N\sum_{k=1}^{n_i}(\alpha^i-\bar{\alpha})\langle\bar{y}^i_{k,t} - x^{i*}_{k,\mu}, \nabla_{x^i_k} f^i_\mu(\bar{\mathbf{y}}_t) - \nabla_{x^i_k} f^i_\mu(\mathbf{x}^*_\mu)\rangle\nonumber\\
&\quad\leq2\sum_{i=1}^N\sum_{k=1}^{n_i}|\alpha^i-\bar{\alpha}|\|\bar{y}^i_{k,t} - x^{i*}_{k,\mu}\|\|\nabla_{x^i_k} f^i_\mu(\bar{\mathbf{y}}_t) - \nabla_{x^i_k} f^i_\mu(\mathbf{x}^*_\mu)\|\nonumber\\
&\quad\leq2L\|\bar{\mathbf{y}}_t-\mathbf{x}^*_\mu\|\sum_{i=1}^N|\alpha^i-\bar{\alpha}|\sum_{k=1}^{n_i}\|\bar{y}^i_{k,t} - x^{i*}_{k,\mu}\|\nonumber\\
&\quad\leq2L\|\bar{\mathbf{y}}_t-\mathbf{x}^*_\mu\|\sum_{i=1}^N\sqrt{n_i}|\alpha^i-\bar{\alpha}|\|\bar{\mathbf{y}}^i_t - \mathbf{x}^{i*}_{\mu}\|\leq2\sqrt{n}L\epsilon_{\alpha}\bar{\alpha}\|\bar{\mathbf{y}}_t-\mathbf{x}^*_\mu\|^2.\label{eq:optimality_gap_expand_4_break_1_result}
\end{align}
For \eqref{eq:optimality_gap_expand_4_break_2},
\begin{align}
&-2\bar{\alpha}\sum_{i=1}^N\sum_{k=1}^{n_i}\langle\bar{y}^i_{k,t} - x^{i*}_{k,\mu}, \nabla_{x^i_k} f^i_\mu(\bar{\mathbf{y}}_t) - \nabla_{x^i_k} f^i_\mu(\mathbf{x}^*_\mu)\rangle\nonumber\\
&\quad=-2\bar{\alpha}\langle \bar{\mathbf{y}}_t - \mathbf{x}^*_\mu, \Phi_\mu(\bar{\mathbf{y}}_t) - \Phi_\mu(\mathbf{x}^*_\mu)\rangle\leq-2\chi\bar{\alpha}\|\bar{\mathbf{y}}_t-\mathbf{x}^*_\mu\|^2.\label{eq:optimality_gap_expand_4_break_2_result}
\end{align}
Combining \eqref{eq:optimality_gap_expand_4_break_1_result} and \eqref{eq:optimality_gap_expand_4_break_2_result}, we obtain for \eqref{eq:optimality_gap_expand_4} that
\begin{align}
&-2\sum_{i=1}^N\sum_{k=1}^{n_i}\alpha^i\langle\bar{y}^i_{k,t} - x^{i*}_{k,\mu}, \nabla_{x^i_k} f^i_\mu(\bar{\mathbf{y}}_t) - \nabla_{x^i_k} f^i_\mu(\mathbf{x}^*_\mu)\rangle\leq2(\sqrt{n}L\epsilon_{\alpha}-\chi)\bar{\alpha}\|\bar{\mathbf{y}}_t-\mathbf{x}^*_\mu\|^2. \label{eq:optimality_gap_expand_4_result}
\end{align}
Finally, taking the conditional expectation for \eqref{eq:optimality_gap_expand} on $\mathcal{H}_t$, and  substituting \eqref{eq:optimality_gap_expand_1_result}, \eqref{eq:optimality_gap_expand_2_result}, \eqref{eq:optimality_gap_expand_3_result} and \eqref{eq:optimality_gap_expand_4_result} into it, we obtain the desired result.
\end{Proof}

Finally, we derive an inequality relation for the third term, $\sum_{i=1}^N\sum_{k=1}^{n_i}\|\varphi^i_{k,t+1} - \mathbf{1}_{n_i}\bar{\varphi}^i_{k,t}\|^2$ in Lemma~\ref{lemma:grad_track_error}.

\begin{Lemma}\label{lemma:grad_track_error}
Under Assumptions~\ref{assumption_graph} and \ref{assumption_local_f_lipschitz}, the total gradient tracking error $\sum_{i=1}^N\sum_{k=1}^{n_i}\|\varphi^i_{k,t} - \mathbf{1}_{n_i}\bar{\varphi}^i_{k,t}\|^2$ satisfies
\begin{align*}
&\sum_{i=1}^N\sum_{k=1}^{n_i}\mathbb{E}[\|\varphi^i_{k,t+1} - \mathbf{1}_{n_i}\bar{\varphi}^i_{k,t+1}\|^2|\mathcal{H}_t]\leq24(n+4)n_s\varsigma L^2\bigg(\frac{3+\bar{\sigma}^2}{2}+\frac{n^2L^2\alpha_{\max}}{\chi}\\
&\quad\quad+24(n+4)n_c\varsigma L^2\alpha_{\max}^2+12n(n+4)L^2\alpha_{\max}^2\bigg)\sum_{i=1}^N\sum_{k=1}^{n_i}\|\mathbf{y}^i_{k,t}-\mathbf{1}_{n_i}\bar{y}^i_{k,t}\|^2\\
&\quad+\bigg(\frac{1+\bar{\sigma}^2}2+48(n+4)n_s\varsigma^2 L^2\alpha_{\max}^2\bigg)\sum_{i=1}^N\sum_{k=1}^{n_i}\mathbb{E}[\|\varphi^i_{k,t} - \mathbf{1}_{n_i}\bar{\varphi}^i_{k,t}\|^2|\mathcal{H}_t]\\
&\quad+24(n+4)n_s\varsigma L^2[2+12n(n+4)L^2\alpha_{\max}^2+24(n+4)n_c\varsigma L^2\alpha_{\max}^2]\|\bar{\mathbf{y}}_t-\mathbf{x}^*_\mu\|^2\\
&\quad+24(n+4)n_s\varsigma L^2[24(n+4)n_c\varsigma G^2+6(n+4)^3n_c\varsigma \mu^2L^2+12n(n+4)G^2\\
&\quad\quad+3n(n+4)^3\mu^2L^2]\alpha_{\max}^2+12(n+4)^3n_s\varsigma \mu^2L^2+48(n+4)n_s\varsigma G^2.
\end{align*}
\end{Lemma}
\begin{Proof}
It is obtained from \eqref{eq:update_grad_track_compact} that
\begin{align*}
\|\varphi^i_{k,t+1} - \mathbf{1}_{n_i}\bar{\varphi}^i_{k,t+1}\|^2&= \|\mathcal{A}^i\varphi^i_{k,t} - \mathbf{1}_{n_i}\bar{\varphi}^i_{k,t}\|^2 + \bigg\|\bigg(I_{n_i}-\frac1{n_i}\mathbf{1}_{n_i}\mathbf{1}_{n_i}^\top\bigg)(\mathbf{g}^i_k(\mathbf{x}_{t+1}) - \mathbf{g}^i_k(\mathbf{x}_t))\bigg\|^2 \\
&\quad+2\bigg\langle \mathcal{A}^i\varphi^i_{k,t} -\mathbf{1}_{n_i}\bar{\varphi}^i_{k,t}, \bigg(I_{n_i}-\frac1{n_i}\mathbf{1}_{n_i}\mathbf{1}_{n_i}^\top\bigg)(\mathbf{g}^i_k(\mathbf{x}_{t+1}) - \mathbf{g}^i_k(\mathbf{x}_t) )\bigg\rangle.
\end{align*}
It is noted that $\|I_{n_i}-\frac1{n_i}\mathbf{1}_{n_i}\mathbf{1}_{n_i}^\top\| = 1$.
Taking the conditional expectation on $\mathcal{H}_t$ yields
\begin{align}
&\mathbb{E}[\|\varphi^i_{k,t+1} - \mathbf{1}_{n_i}\bar{\varphi}^i_{k,t+1}\|^2|\mathcal{H}_t] \leq \sigma_{\mathcal{A}^i}^2\mathbb{E}[\|\varphi^i_{k,t} - \mathbf{1}_{n_i}\bar{\varphi}^i_{k,t}\|^2|\mathcal{H}_t] + \mathbb{E}[\|\mathbf{g}^i_k(\mathbf{x}_{t+1}) - \mathbf{g}^i_k(\mathbf{x}_t)\|^2|\mathcal{H}_t]\nonumber\\
&\quad\quad+2\mathbb{E}[\| \mathcal{A}^i\varphi^i_{k,t} - \mathbf{1}_{n_i}\bar{\varphi}^i_{k,t}\| \|\mathbf{g}^i_k(\mathbf{x}_{t+1}) - \mathbf{g}^i_k(\mathbf{x}_t)\||\mathcal{H}_t]\nonumber\\
&\quad\leq \sigma_{\mathcal{A}^i}^2\mathbb{E}[\|\varphi^i_{k,t} - \mathbf{1}_{n_i}\bar{\varphi}^i_{k,t}\|^2|\mathcal{H}_t] + \mathbb{E}[\|\mathbf{g}^i_k(\mathbf{x}_{t+1}) - \mathbf{g}^i_k(\mathbf{x}_t)\|^2|\mathcal{H}_t]\nonumber\\
&\quad\quad +\frac{1-\sigma_{\mathcal{A}^i}^2}2\mathbb{E}[\|\varphi^i_{k,t} - \mathbf{1}_{n_i}\bar{\varphi}^i_{k,t}\|^2|\mathcal{H}_t]+\frac{2\sigma_{\mathcal{A}^i}^2}{1-\sigma_{\mathcal{A}^i}^2}\mathbb{E}[\|\mathbf{g}^i_k(\mathbf{x}_{t+1}) - \mathbf{g}^i_k(\mathbf{x}_t)\|^2|\mathcal{H}_t]\nonumber\\
&\quad\leq\frac{1+\bar{\sigma}^2}2\mathbb{E}[\|\varphi^i_{k,t} - \mathbf{1}_{n_i}\bar{\varphi}^i_{k,t}\|^2|\mathcal{H}_t]+\varsigma\mathbb{E}[\|\mathbf{g}^i_k(\mathbf{x}_{t+1}) - \mathbf{g}^i_k(\mathbf{x}_t)\|^2|\mathcal{H}_t].\label{eq:grad_track_error_term123}
\end{align}
The last term of \eqref{eq:grad_track_error_term123} follows from \eqref{eq:pi_i_jk_bound} that
\begin{align*}
&\mathbb{E}[\|\mathbf{g}^i_k(\mathbf{x}_{t+1}) - \mathbf{g}^i_k(\mathbf{x}_t)\|^2|\mathcal{H}_t]\leq 2\sum_{j=1}^{n_i}(\mathbb{E}[\|g^i_{jk}(\mathbf{x}_{t+1})\|^2|\mathcal{H}_t]+\mathbb{E}[\|g^i_{jk}(\mathbf{x}_t)\|^2|\mathcal{H}_t])\\
&\quad\leq 24n_i(n+4)L^2\sum_{i=1}^N\sum_{k=1}^{n_i}\mathbb{E}[\|\mathbf{y}^i_{k,t+1}-\mathbf{1}_{n_i}\bar{y}^i_{k,t+1}\|^2|\mathcal{H}_t]+24n_i(n+4)L^2\mathbb{E}[\|\bar{\mathbf{y}}_{t+1}-\mathbf{x}^*_\mu\|^2|\mathcal{H}_t]\\
&\quad\quad+24n_i(n+4)L^2\sum_{i=1}^N\sum_{k=1}^{n_i}\|\mathbf{y}^i_{k,t}-\mathbf{1}_{n_i}\bar{y}^i_{k,t}\|^2+24n_i(n+4)L^2\|\bar{\mathbf{y}}_t-\mathbf{x}^*_\mu\|^2+12n_i(n+4)^3\mu^2L^2\\
&\quad\quad+48n_i(n+4)G^2.
\end{align*}
Invoking Lemmas~\ref{lemma:conensus} and \ref{lemma:optimality_gap} in the above relation, and summing \eqref{eq:grad_track_error_term123} over $k\in\mathcal{V}^i,i\in\mathcal{N}$ complete the proof.
\end{Proof}

\subsection{Proof of Theorem~\ref{theorem:optimality}}\label{subsec:proof_of_theorem}

Now, we proceed to the proof of Theorem~\ref{theorem:optimality}. Based on the results in Lemmas~\ref{lemma:conensus}, \ref{lemma:optimality_gap} and \ref{lemma:grad_track_error}, we can construct a linear system by taking the total expectation on the corresponding relations.
\begin{align}
\Psi_{t+1} \leq \mathbf{M}\Psi_t + \Upsilon, \label{eq:dynamical_system}
\end{align}
where
\begin{align*}
\Psi_t &\triangleq \begin{bmatrix} \sum_{i=1}^N\sum_{k=1}^{n_i}\mathbb{E}[\|\mathbf{y}^i_{k,t}-\mathbf{1}_{n_i}\bar{y}^i_{k,t}\|_{\bm{\nu}^i_r}^2]\\\mathbb{E}[\|\bar{\mathbf{y}}_t-\mathbf{x}^*_\mu\|^2]\\\sum_{i=1}^N\sum_{k=1}^{n_i}\mathbb{E}[\|\varphi^i_{k,t} - \mathbf{1}_{n_i}\bar{\varphi}^i_{k,t}\|^2] \end{bmatrix},\Upsilon \triangleq \begin{bmatrix}m_{12}{\alpha}_{\max}^2\\ m_{13}{\alpha}_{\max}^2 \\ m_{14}+m_{15}{\alpha}_{\max}^2\end{bmatrix},\\
\mathbf{M} &\triangleq\begin{bmatrix}1-m_1+m_2{\alpha}_{\max}^2 & m_2{\alpha}_{\max}^2 & m_3{\alpha}_{\max}^2\\ m_4{\alpha}_{\max}+m_5{\alpha}_{\max}^2 & 1-m_6\bar{\alpha}+m_5{\alpha}_{\max}^2 & 0 \\ m_7+m_8{\alpha}_{\max}+m_9{\alpha}_{\max}^2 & m_{10}+m_9{\alpha}_{\max}^2 & 1-m_1+m_{11}{\alpha}_{\max}^2\end{bmatrix},
\end{align*}
$m_1\triangleq\frac{1-\bar{\sigma}^2}{2}$, $m_2\triangleq 24(n+4)n_c\varsigma L^2$, $m_3\triangleq 2\varsigma$, $m_4 \triangleq n^2L^2/\chi$, $m_5\triangleq 12n(n+4) L^2$, $m_6 \triangleq \chi-2\sqrt{n}L\epsilon_{\alpha}$, $m_7\triangleq 12(n+4)n_s\varsigma L^2(3+\bar{\sigma}^2)$, $m_8\triangleq 24(n+4)n_s\varsigma L^2m_4$, $m_9\triangleq 24(n+4)n_s\varsigma L^2(m_2+m_5)$, $m_{10}\triangleq 48(n+4)n_s\varsigma L^2$, $m_{11}\triangleq \varsigma m_{10}$, $m_{12} \triangleq 24(n+4)n_c\varsigma G^2+6(n+4)^3n_c\varsigma \mu^2L^2$, $m_{13} \triangleq 12n(n+4)G^2+3n(n+4)^3 \mu^2L^2$, $m_{14}\triangleq 48(n+4)n_s\varsigma G^2+12(n+4)^3n_s\varsigma \mu^2L^2$ and $m_{15}\triangleq 24(n+4)n_s\varsigma L^2(m_{12}+m_{13})$.


For the linear system \eqref{eq:dynamical_system}, we aim to prove $\rho(\mathbf{M})<1$ such that each component of $\Psi_t$ can linearly converge to a neighborhood of 0 \cite{Horn1990}. 

We adopt the following result to guarantee $\rho(\mathbf{M})<1$:
\begin{Lemma}\label{lemma:matrix_spectral_radius}
(see \cite[Cor.~8.1.29]{Horn1990}) Let $A\in\mathbb{R}^{m\times m}$ be a matrix with non-negative entries and $\bm{\nu}\in\mathbb{R}^m$ be a vector with positive entries. If there exists a constant $\lambda\geq0$ such that $A\bm{\nu} < \lambda \bm{\nu}$, then $\rho(A) < \lambda$.
\end{Lemma}

To apply Lemma~\ref{lemma:matrix_spectral_radius}, each element of $\mathbf{M}$ should be non-negative. Hence, we may set $m_6>0$ and ${\alpha}_{\max}<\frac1{m_6}$, \textit{i.e.},
\begin{align*}
{\alpha}_{\max}<\frac1{m_6}, \epsilon_{\alpha} < \frac{\chi}{2\sqrt{n}L}.
\end{align*}
Next, based on Lemma~\ref{lemma:matrix_spectral_radius}, it suffices to find a vector $\bm{\nu} \triangleq [\nu_1,\nu_2,\nu_3]^\top$ with $\nu_1,\nu_2,\nu_3>0$ such that $\mathbf{M}_\alpha\bm{\nu} < \bm{\nu}$, \textit{i.e.},
\begin{align*}
&(1-m_1+m_2{\alpha}_{\max}^2)\nu_1+(m_2{\alpha}_{\max}^2)\nu_2 + (m_3{\alpha}_{\max}^2)\nu_3 < \nu_1,\\
&(m_4{\alpha}_{\max}+m_5{\alpha}_{\max}^2)\nu_1+ (1-m_6\bar{\alpha}+m_5{\alpha}_{\max}^2)\nu_2 < \nu_2,\\
&(m_7+m_8{\alpha}_{\max}+m_9{\alpha}_{\max}^2)\nu_1 + (m_{10}+m_9{\alpha}_{\max}^2)\nu_2 + (1-m_1+m_{11}{\alpha}_{\max}^2)\nu_3 < \nu_3.
\end{align*}
Without loss of generality, we may set $\nu_3 = 1$. It remains to find $\nu_1$ and $\nu_2$ such that the following inequalities hold
\begin{subequations}\label{eq:ineq_t}
\begin{align}
&(m_2\nu_1+m_2\nu_2+m_3){\alpha}_{\max}^2 < m_1\nu_1, \label{eq:ineq_t1}\\
&(m_5\nu_1+m_5\nu_2){\alpha}_{\max} < \frac{m_6\nu_2}n-m_4\nu_1, \label{eq:ineq_t2}\\
&(m_9\nu_1+m_9\nu_2+m_{11}){\alpha}_{\max}^2<m_1-(m_7+m_8)\nu_1-m_{10}\nu_2,\label{eq:ineq_t3}
\end{align}
\end{subequations}
where we have applied $\frac{\bar{\alpha}}{{\alpha}_{\max}}>\frac1n$ in \eqref{eq:ineq_t2}, and forced ${\alpha}_{\max}<1$ in \eqref{eq:ineq_t3}.

To ensure the existence of ${\alpha}_{\max}$, the RHS of \eqref{eq:ineq_t} has to be positive. Hence, we may set 
\begin{align*}
\nu_1=\frac{m_1m_6}{4nm_4m_{10}+2m_6m_7+2m_6m_8},\nu_2=\frac{nm_1m_4}{2nm_4m_{10}+m_6m_7+m_6m_8}. 
\end{align*}
Then, the three inequalities in \eqref{eq:ineq_t} can be solved, which gives
\begin{align*}
&{\alpha}_{\max}<\alpha_1,{\alpha}_{\max}<\alpha_2,{\alpha}_{\max}<\alpha_3,
\end{align*}
where 
\begin{align*}
&\alpha_1\triangleq \sqrt{\frac{m_1^2m_6}{m_1m_2m_6+2nm_1m_2m_4+4nm_3m_4m_{10}+2m_3m_6m_7+2m_3m_6m_8}},\\
&\alpha_2\triangleq \frac{m_1m_4m_6}{m_1m_5m_6+2nm_1m_4m_5},\\
&\alpha_3\triangleq \sqrt{\frac{m_1(2nm_4m_{10}+m_6m_7+m_6m_8)}{m_1m_6m_9+2nm_1m_4m_9+m_{11}(4nm_4m_{10}+2m_6m_7+2m_6m_8)}}.
\end{align*}
Therefore, the range of the step-size is given by
\begin{align*}
0<{\alpha}_{\max}<\min\bigg\{\alpha_1,\alpha_2,\alpha_3,\frac1{m_6},1\bigg\}, \epsilon_{\alpha} < \frac{\chi}{2\sqrt{n}L}.
\end{align*}

Furthermore, taking the limsup on both sides of \eqref{eq:dynamical_system}
\begin{align*}
\limsup_{t\to\infty}\Psi_t \leq \mathbf{M}\limsup_{t\to\infty}\Psi_t + \Upsilon,
\end{align*}
which gives
\begin{align*}
(I_3 - \mathbf{M})\limsup_{t\to\infty}\Psi_t \leq \Upsilon,
\end{align*}
where
\begin{align*}
I_3 - \mathbf{M} = \begin{bmatrix}m_1-m_2{\alpha}_{\max}^2 & -m_2{\alpha}_{\max}^2 & -m_3{\alpha}_{\max}^2\\ -m_4{\alpha}_{\max}-m_5{\alpha}_{\max}^2 & m_6\bar{\alpha}-m_5{\alpha}_{\max}^2 & 0 \\ -m_7-m_8{\alpha}_{\max}-m_9{\alpha}_{\max}^2 & -m_{10}-m_9{\alpha}_{\max}^2 & m_1-m_{11}{\alpha}_{\max}^2\end{bmatrix}.
\end{align*}
It can be obtained that
\begin{align*}
det(I_3 - \mathbf{M}) &\triangleq (m_1-m_{11}{\alpha}_{\max}^2)[(m_1-m_2{\alpha}_{\max}^2)(m_6\bar{\alpha}-m_5{\alpha}_{\max}^2)\\
&\quad-m_2{\alpha}_{\max}^2(m_4{\alpha}_{\max}+m_5{\alpha}_{\max}^2)]\\
&>(m_1-m_{11}{\alpha}_{\max}^2)\bigg[(m_1-m_2{\alpha}_{\max}^2)\bigg(\frac{m_6{\alpha}_{\max}}n-m_5{\alpha}_{\max}^2\bigg)\\
&\quad-m_2{\alpha}_{\max}^2(m_4{\alpha}_{\max}+m_5{\alpha}_{\max}^2)\bigg]\\
&={\alpha}_{\max}(m_1-m_2{\alpha}_{\max}^2)\bigg[\frac{m_1m_6}n-m_1m_5{\alpha}_{\max}-m_2\bigg(m_4+\frac{m_6}n\bigg){\alpha}_{\max}^2\bigg],\\
adj(I_3 - \mathbf{M})_{11} &\triangleq (m_1-m_{11}{\alpha}_{\max}^2)(m_6\bar{\alpha}-m_5{\alpha}_{\max}^2)\leq{\alpha}_{\max}(m_1-m_{11}{\alpha}_{\max}^2)(m_6-m_5{\alpha}_{\max}),\\
adj(I_3 - \mathbf{M})_{12} &\triangleq {\alpha}_{\max}^2[m_2(m_1-m_{11}{\alpha}_{\max}^2)+m_3(m_{10}+m_9{\alpha}_{\max}^2)],\\
adj(I_3 - \mathbf{M})_{13} &\triangleq m_3{\alpha}_{\max}^2(m_6\bar{\alpha}-m_5{\alpha}_{\max}^2)\leq m_3{\alpha}_{\max}^3(m_6-m_5{\alpha}_{\max}),\\
adj(I_3 - \mathbf{M})_{21} &\triangleq {\alpha}_{\max}(m_1-m_{11}{\alpha}_{\max}^2)(m_4+m_5{\alpha}_{\max}),\\
adj(I_3 - \mathbf{M})_{22} &\triangleq (m_1-m_2{\alpha}_{\max}^2)(m_1-m_{11}{\alpha}_{\max}^2)-m_3{\alpha}_{\max}^2(m_7+m_8{\alpha}_{\max}+m_9{\alpha}_{\max}^2),\\
adj(I_3 - \mathbf{M})_{23} &\triangleq m_3{\alpha}_{\max}^3(m_4+m_5{\alpha}_{\max})
\end{align*}
Then, we have
\begin{align*}
&\limsup_{t\to\infty}\sum_{i=1}^N\sum_{k=1}^{n_i}\mathbb{E}[\|\mathbf{y}^i_{k,t}-\mathbf{1}_{n_i}\bar{y}^i_{k,t}\|^2]\leq[(I_3 - \mathbf{M})^{-1}\Upsilon]_1\\
&\quad=\frac{adj(I_3 - \mathbf{M})_{11}[\Upsilon]_1}{det(I_3 - \mathbf{M})}+\frac{adj(I_3 - \mathbf{M})_{12}[\Upsilon]_2}{det(I_3 - \mathbf{M})}+\frac{adj(I_3 - \mathbf{M})_{13}[\Upsilon]_3}{det(I_3 - \mathbf{M})}=\mathcal{O}({\alpha}_{\max}^2),
\end{align*}
and
\begin{align*}
&\limsup_{t\to\infty}\mathbb{E}[\|\bar{\mathbf{y}}_t-\mathbf{x}^*_\mu\|^2]\leq[(I_3 - \mathbf{M})^{-1}\Upsilon]_2\\
&\quad=\frac{adj(I_3 - \mathbf{M})_{21}[\Upsilon]_1}{det(I_3 - \mathbf{M})}+\frac{adj(I_3 - \mathbf{M})_{22}[\Upsilon]_2}{det(I_3 - \mathbf{M})}+\frac{adj(I_3 - \mathbf{M})_{23}[\Upsilon]_3}{det(I_3 - \mathbf{M})}=\mathcal{O}({\alpha}_{\max}).
\end{align*}
Thus, it follows from \eqref{eq:x_minus_ybar} that
\begin{align*}
&\limsup_{t\to\infty}\mathbb{E}[\|\mathbf{x}_t - \mathbf{x}^*_\mu\|^2]\leq2\limsup_{t\to\infty}\mathbb{E}[\|\mathbf{x}_t - \bar{\mathbf{y}}_t\|^2]+2\limsup_{t\to\infty}\mathbb{E}[\|\bar{\mathbf{y}}_t - \mathbf{x}^*_\mu\|^2]\\
&\quad\leq2\limsup_{t\to\infty}\sum_{i=1}^N\sum_{k=1}^{n_i}\mathbb{E}[\|\mathbf{y}^i_{k,t}-\mathbf{1}_{n_i}\bar{y}^i_{k,t}\|^2]+2\limsup_{t\to\infty}\mathbb{E}[\|\bar{\mathbf{y}}_t - \mathbf{x}^*_\mu\|^2]=\mathcal{O}({\alpha}_{\max}).
\end{align*}
Invoking Lemma~\ref{lemma:NE_gap} yields
\begin{align*}
&\limsup_{t\to\infty}\mathbb{E}[\|\mathbf{x}_t - \mathbf{x}^*\|^2]\leq2\limsup_{t\to\infty}\mathbb{E}[\|\mathbf{x}_t - \mathbf{x}^*_\mu\|^2]+2\|\mathbf{x}^*_\mu- \mathbf{x}^*\|^2=\mathcal{O}({\alpha}_{\max})+\mathcal{O}(\mu),
\end{align*}
which completes the proof.

\section{Numerical Simulations}\label{sec:simulation}

We illustrate the proposed NE seeking strategy on a connectivity control game \cite{Stankovic2012}, played among a number of sensor networks.
Specifically, there are $N$ sensor networks, where each sensor network contains $n_i$ sensors. Let $x^i_j=[x^i_{j,1},x^i_{j,2}]^\top\in\mathbb{R}^2$ denote the position of sensor $j$ (referred to as an agent) from a sensor network $i$ (referred to as a cluster). Then, this sensor aims to seek a tradeoff between a local cost, $l^i_j(\mathbf{x}^i)$ (\textit{e.g.}, source seeking and positioning) and the global cost, $h^i_j(\mathbf{x})$ (\textit{e.g.}, connectivity preservation with other sensor networks).
Hence, the cost function to be minimized by this sensor is given by
\begin{equation*}
f^i_j(\mathbf{x}) = l^i_j(\mathbf{x}^i) + h^i_j(\mathbf{x}),
\end{equation*}
where
\begin{align*}
l^i_j(\mathbf{x}^i) &= \mathbf{x}^{i\top}a^i_j\mathbf{x}^i+b^{i\top}_j\mathbf{x}^i+c^i_j,\\
h^i_j(\mathbf{x}) &= \sum_{k\in\mathcal{N}_i}d^i_j\|e^i_{jk}x^i_j-\mathbf{x}^k\|^2,
\end{align*}
and $a^i_j,b^i_j,c^i_j,d^i_j,e^i_{jk}$ are constant matrices or vectors of appropriate dimensions, and $\mathcal{N}_i$ stands for the set of neighbors of sensor network $i$ in a connected graph characterizing their position dependence. Specifically, if $k\in\mathcal{N}_i$, then the corresponding term $\|e^i_{jk}x^i_j-\mathbf{x}^k\|^2$ represents the intention of sensor $j$ from a sensor network $i$ to preserve the connectivity with the sensors from sensor network $k$.

\begin{figure}[!t]
\centering
\includegraphics[width=2.6in]{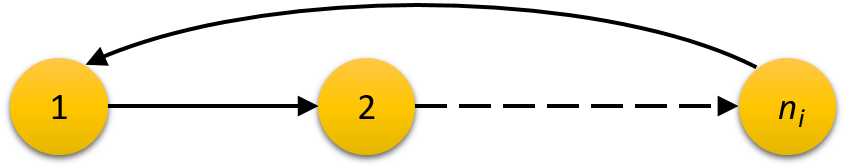}  %
\caption{Communication network.}
\label{fig:network.png}%
\end{figure} 

\begin{figure*}[!htb]
\centering
\subfloat[Network 1]{\includegraphics[width=2.1in]{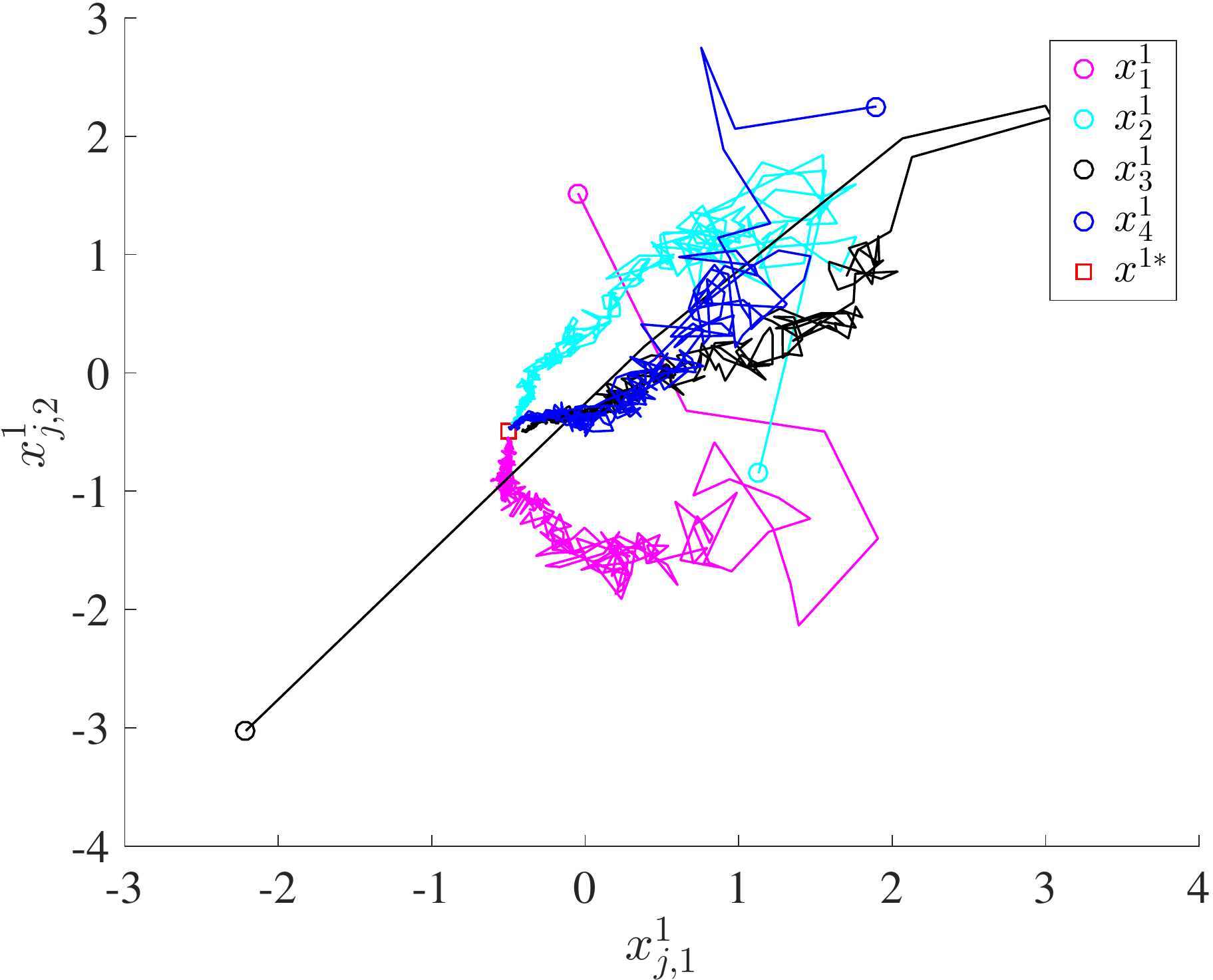}}
\hfil
\subfloat[Network 2]{\includegraphics[width=2.1in]{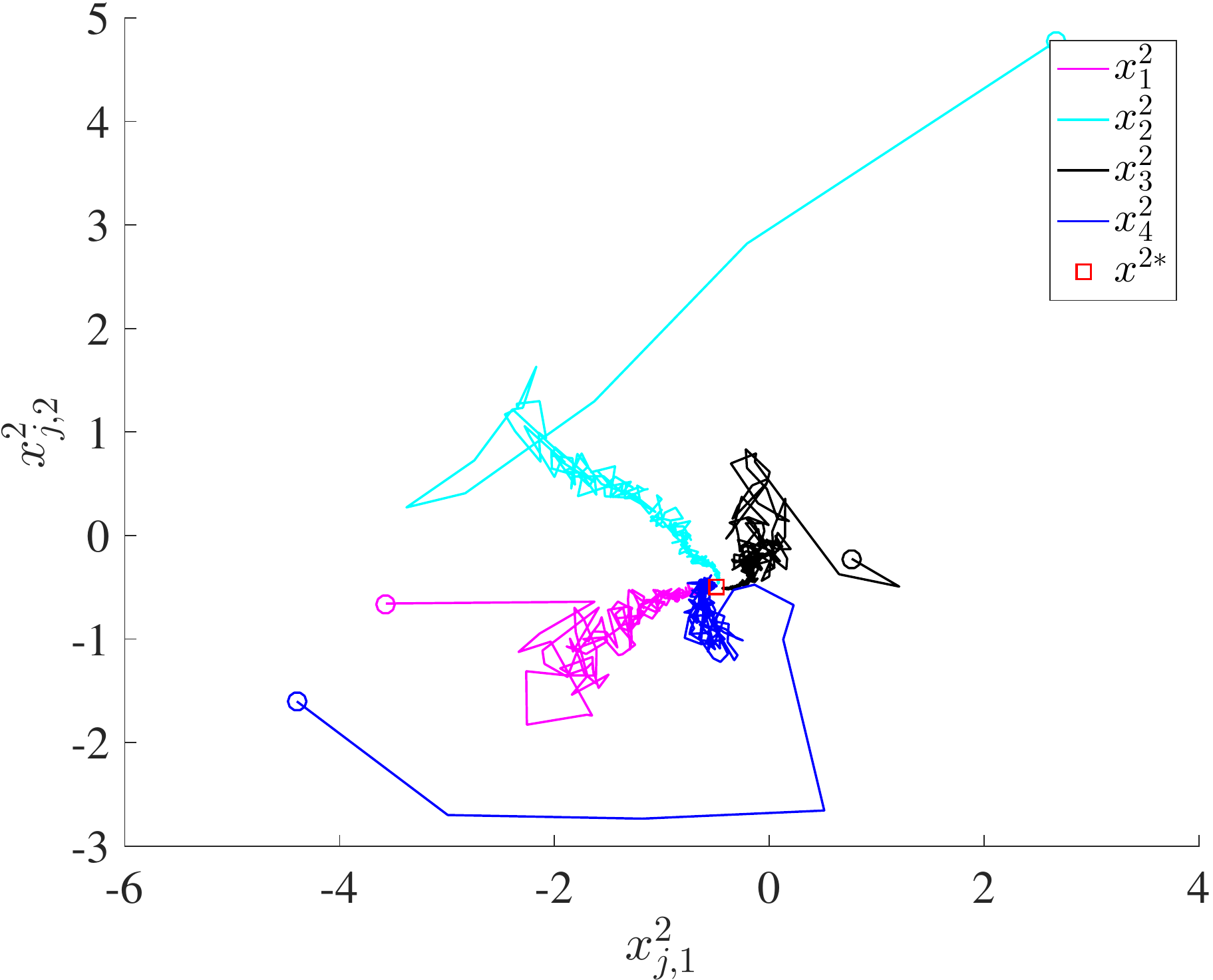}}
\hfil
\subfloat[Network 3]{\includegraphics[width=2.1in]{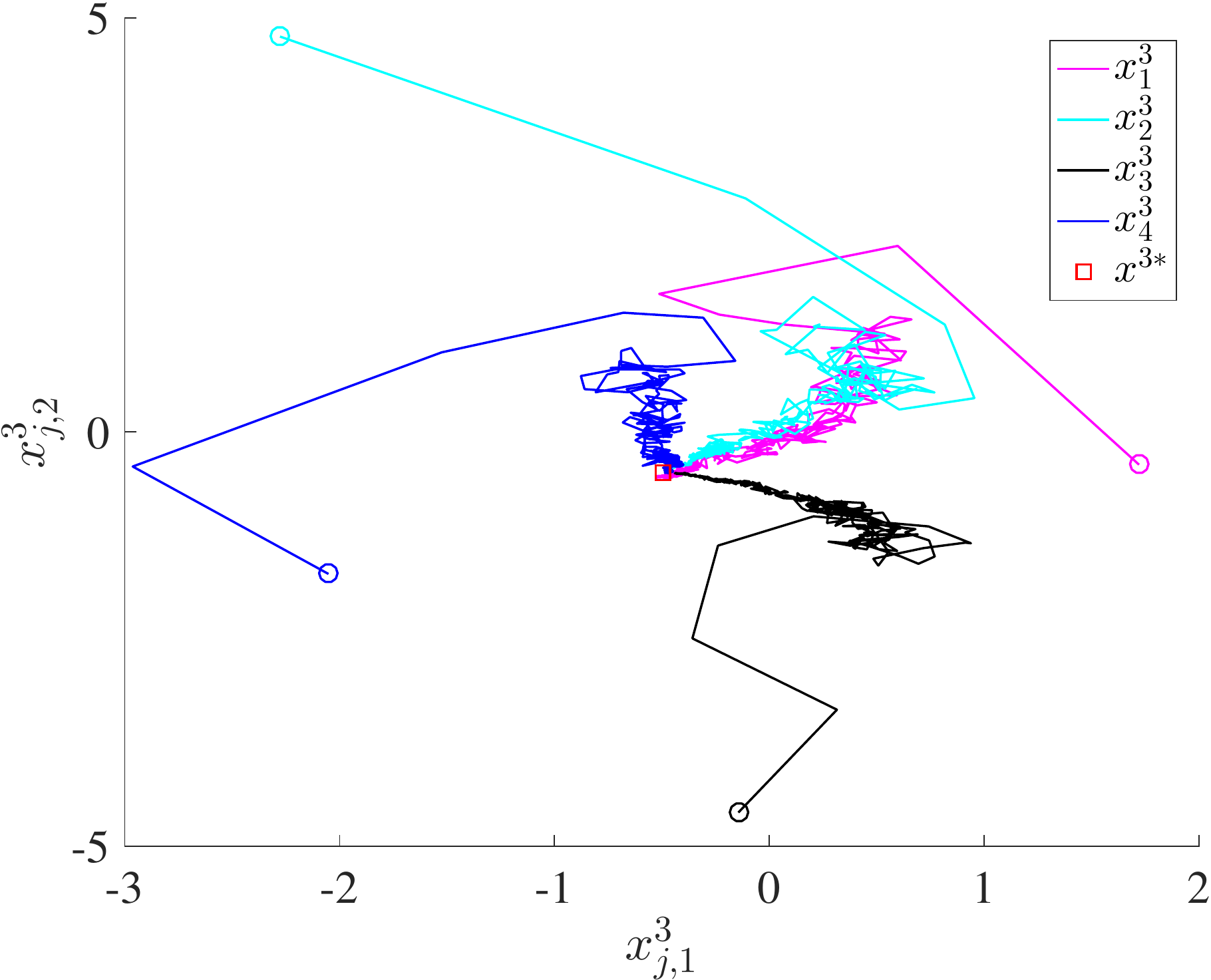}}
\caption{Trajectories of sensors' positions in different networks.}
\label{fig:action_rgf_con_ex2.png}
\end{figure*}

In this simulation, we consider $N = 3$ and $n_i = 4$. The local and global costs are set as $l^i_j = i[\|x^i_j\|^2+\mathbf{1}^\top_2x^i_j+j]$ for $j=1,\ldots,4$, $i=1,2,3$ and $h^1_j = \|x^1_j-x^2_j\|^2$, $h^2_j = \|x^2_j-x^3_j\|^2$, $h^3_j = \|x^3_j-x^1_j\|^2$ for $j=1,\ldots,4$. Then, it is readily verified that Assumptions~\ref{assumption_local_f_lipschitz} and \ref{assumption_game_mapping} hold. The directed communication graph for each sensor network $i$ is as shown in Fig.~\ref{fig:network.png}.
For the algorithm parameters, we let the smoothing parameter be $\mu = 10^{-4}$, and the constant step-sizes for sensors of network $i$ be $\alpha^i = 0.1$, $0.08$, $0.06$, respectively. Thus, $\epsilon_\alpha = 0.2041$. We initialize the algorithm with arbitrary $x^i_{j,0}$, $y^i_{jk,0}$ and $\varphi^i_{jk,0} = {g}^i_{jk}(\mathbf{x}_0)$. The trajectories of the sensors' positions for the three sensor networks are plotted in Fig.~\ref{fig:action_rgf_con_ex2.png}. It can be seen that the positions of all sensors can almost converge to the NE. Also, more `zigzags' can be observed for the case of a larger step-size, since the update is more aggressive.

\begin{figure*}[!htb]
\centering
\subfloat[Fixed heterogeneity]{\includegraphics[width=2.6in]{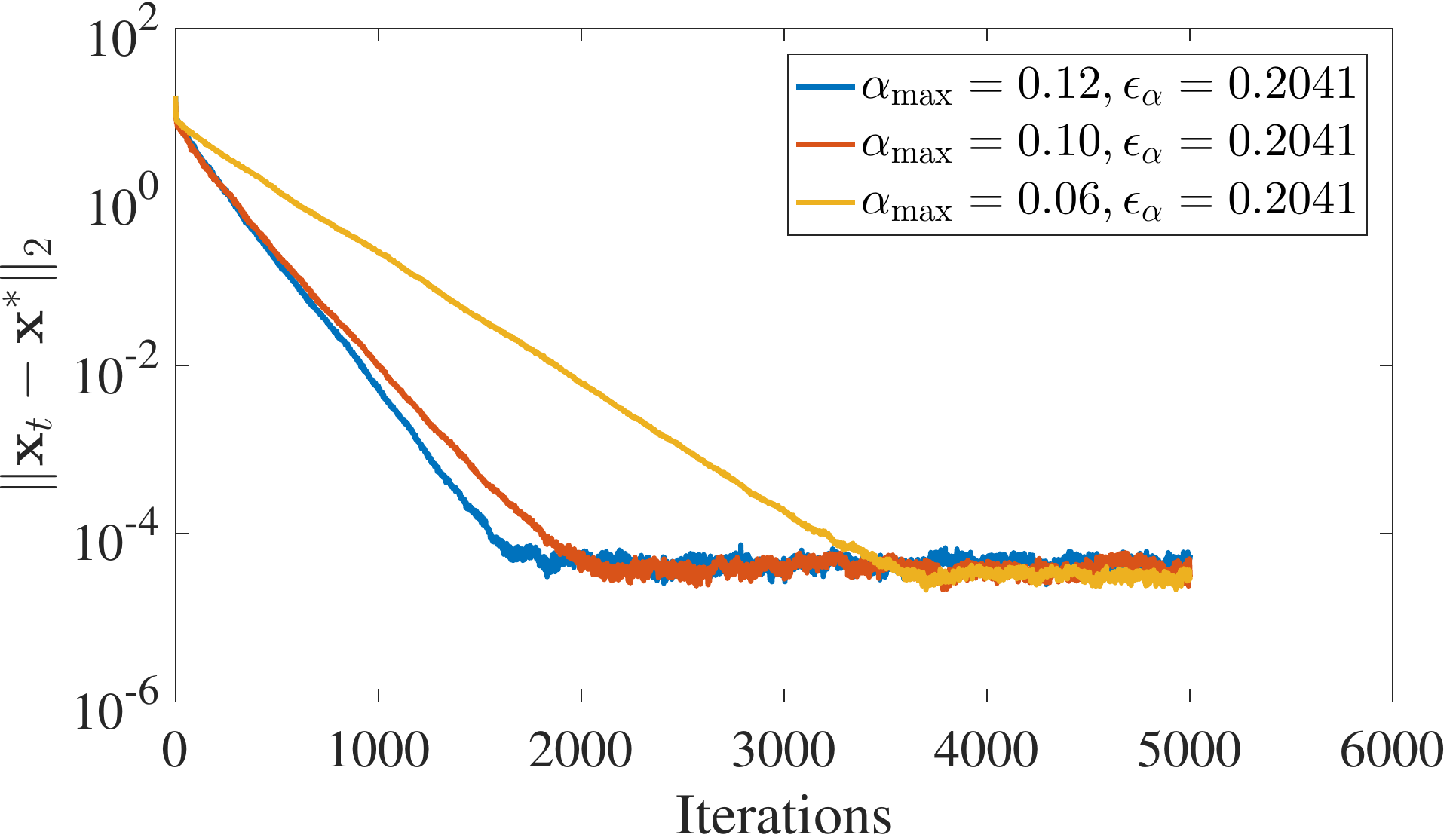}
\label{fig:convergence_rate_alpha_max_ex2.png}}
\hfil
\subfloat[Fixed largest step-size]{\includegraphics[width=2.6in]{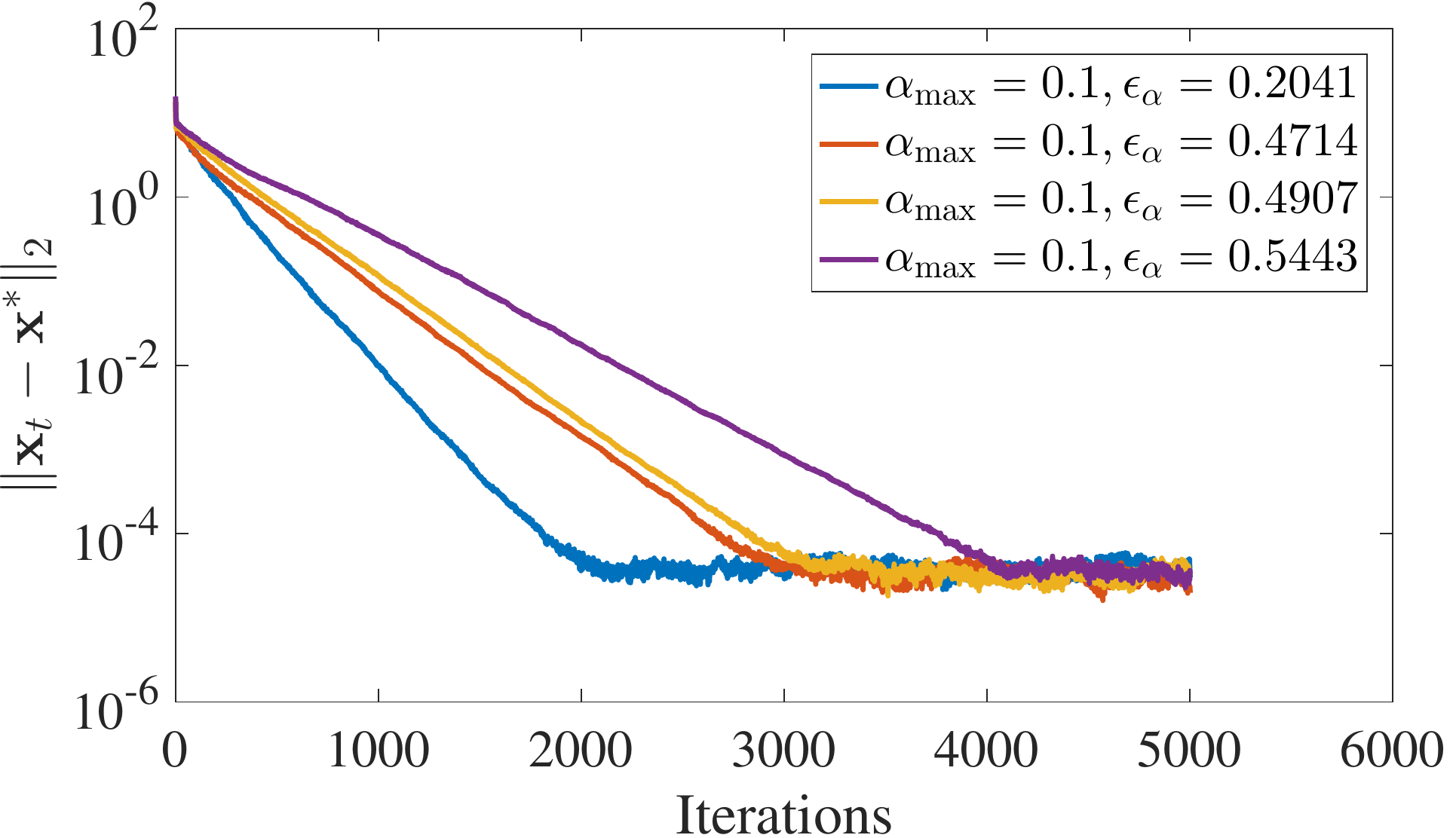}
\label{fig:convergence_rate_ratio_ex2.png}}
\caption{Trajectories of the error gap $\|\mathbf{x}_t - \mathbf{x}^*\|$.}
\end{figure*}

Next, we illustrate the convergence rate results. First, we set the constant step-sizes for sensors of network $i$ be $\alpha^i_j = 0.1a$, $0.08a$ and $0.06a$, respectively, and let $a = 1.2$, $1$ and $0.6$, respectively. Hence, we  fix the heterogeneity of the step-size $\epsilon_\alpha = 0.2041$, and set the largest step-size to ${\alpha}_{\max} = 0.12$, $0.1$ and $0.06$, respectively. The trajectories of the error gap $\|\mathbf{x}_t - \mathbf{x}^*\|$ with these settings are plotted in Fig.~\ref{fig:convergence_rate_alpha_max_ex2.png}. Then, we fix the largest step-size to ${\alpha}_{\max} = 0.1$ and the averaged step-size $\bar{\alpha} = 0.06$, and set the heterogeneity of the step-size $\epsilon_\alpha = 0.2041$, $0.4714$, $0.4907$ and $0.5443$, respectively. The trajectories of the error gap $\|\mathbf{x}_t - \mathbf{x}^*\|$ with these settings are plotted in Fig.~\ref{fig:convergence_rate_ratio_ex2.png}. As can be seen from both figures, the error gap descends linearly for all cases. Moreover, the convergence speed is faster with larger step-sizes and smaller heterogeneity, which verifies the derived results in Theorem~\ref{theorem:optimality}. 

\section{Conclusions}\label{sec:conclusion}
This work has studied an $N$-cluster non-cooperative game problem, where the agents' cost functions are possibly non-smooth and the explicit expressions are unknown. By integrating the Gaussian smoothing techniques with the gradient tracking, a gradient-free NE seeking algorithm has been developed, in which the agents are allowed to select their own preferred constant step-sizes. We have shown that, when the largest step-size is sufficiently small, the agents' actions approximately converge to the unique NE under a strongly monotone game mapping condition, and the error gap is proportional to the largest step-size and the smoothing parameter. 
Finally, the derived results have been verified by numerical simulations.








\bibliographystyle{IEEEtran}
\bibliography{d_ne_coalition_rgf_reference}

\end{document}